\documentclass{article}
\usepackage[utf8]{inputenc}
\usepackage[usenames,dvipsnames]{xcolor}
\usepackage{amsmath,amssymb,graphicx,tikz,todonotes}
\usepackage{hyperref}
\usepackage[margin=1in]{geometry}
\usepackage{layouts}
\usepackage[titletoc,title]{appendix} 
\usepackage{tabularx}
\usepackage{multirow}
\usepackage{placeins}
\usepackage{mathtools}
\usepackage[normalem]{ulem}
\usepackage{xcolor}
\usepackage{relsize}
\usepackage[style=numeric, citestyle=numeric, backend=biber]{biblatex}

\usepackage{verbatim}

\usepackage[font=small,labelfont=bf]{caption}
\newcommand{\bi}[1]{\boldsymbol{#1}}

\def\dV0{{{\rm d}V_0}}
\def\V0{{V_0}}
\def\Vi{{V_{0,i}}}
\def\dA0{{\mathrm{d} A_0}}
\def\A0{{\partial V_0}}
\def\Ai{{\partial V_{0,i}}}

\def\sumalpha{\sum_{\alpha=1}^{n_\alpha}}
\def\sumbeta{\sum_{\beta=1}^{n_\alpha}}
\def\setka{\{k_\alpha\}_{\alpha=1}^{n_\alpha}}
\def\setnablaka{\{\boldsymbol{\nabla}_0 k_\alpha\}_{\alpha=1}^{n_\alpha}}

\def\sumgrains{\mathlarger{\sum_{i=1}^{n_\mathrm{grains}}}}

\title{Phase-field modeling of ductile fracture across grain boundaries in polycrystals }

\author{
Kim Louisa Auth\textsuperscript{a,}\footnote{Corresponding author: \texttt{kim.auth@chalmers.se} (Kim Louisa Auth)}, 
Jim Brouzoulis\textsuperscript{b}, Magnus Ekh\textsuperscript{a}}
\date{}

\begin{document}
\maketitle
\noindent
\textsuperscript{a} \textit{\small Division of Material and Computational Mechanics, Department of Industrial and Materials Science, Chalmers University of Technology, 41296 Gothenburg, Sweden}
\newline
\textsuperscript{b} \textit{\small Division of Dynamics, Department of Mechanics and Maritime Sciences, Chalmers University of Technology, 41296 Gothenburg, Sweden}

\section*{Abstract}
In this study, we address damage initiation and micro-crack formation in ductile failure of polycrystalline metals. We show how our recently published thermodynamic framework for ductile phase-field fracture of single crystals can be extended to polycyrstalline structures. A key feature of this framework is that is accounts for size effects by adopting gradient-enhanced (crystal) plasticity. Gradient-enhanced plasticity requires the definition of boundary conditions representing the plastic slip transmission resistance of the boundaries.
In this work, we propose a novel type of micro-flexible boundary condition for gradient-plasticity, which couples the slip transmission resistance with the phase-field damage such that the resistance locally changes during the fracturing process. The formulation permits to maintain the effect of grain boundaries as obstacles for plastic slip during plastification, while also accounting for weakening of their resistance during the softening phase.
In numerical experiments, the new damage-dependent boundary condition is compared to classical micro-free and micro-hard boundary conditions in polycrystals and it is demonstrated that it indeed produces a response that transitions from micro-hard to micro-free as the material fails. We show that the formulation maintains resistance to slip transmission during hardening, but can generate micro-cracks across grain boundaries during the fracture process.
We further show examples of how the model can be used to simulate void coalescence and three-dimensional crack fronts in polycrystals.

\section{Introduction}
Understanding the phenomena underlying ductile fracture is important for efficient design of many engineering components. Ductile fracture in metals usually occurs after significant plastic deformation. The fracture process proceeds from plastification to damage initiation, followed by the development of microscopic cracks and eventually the evolution of those into macroscopic cracks that can cause component failure. 
In order to simulate this process, it is important to develop models and computational tools that are able to include the physical mechanisms on the grain scale.

Gradient-extended crystal plasticity models have for a long time been used to model size-dependent behavior of polycrystals see e.g. \cite{Gurtin2002}, \cite{Evers2004}, \cite{Ekh2007}, \cite{Svendsen2010}. The size-dependence originates from geometrical necessary dislocations at grain boundaries and how they transmit as well as interact across interfaces like grain boundaries which influence crack initiation and propagation, cf. e.g. \cite{Kacher2014}.  In the context of gradient-extended crystal plasticity modeling, slip transmission mechanisms are accounted for by grain boundary conditions. Such conditions have been discussed in the literature for various gradient-extended plasticity models
by introduction of interface energies, see e.g. \cite{Aifantis2005}, \cite{Fredriksson2005}, or by accounting for slip direction mismatch and grain boundary orientation, see e.g. \cite{Gurtin2008}, \cite{Ekh2011}, \cite{Yalcinkaya2021}, \cite{Wulfinghoff2013}, \cite{McBride2016} and the review by \citeauthor{Bayerschen2015} \cite{Bayerschen2015}.

Damage initiates by void nucleation, often at inclusions and/or in shear bands. For ductile fracture, void nucleation frequently occurs inside single grains, even though it can occur at grain boundaries as well \cite{Noell2017}. During the fracture process these voids grow and coalesce, whereby micro-cracks crossing grain boundaries are formed. The arrangement of voids as well as the way they coalesce strongly depends on the type of loading \cite{Pineau2016}, \cite{Noell2018}, \cite{Azman2022}.
Damage and fracture within the grains may occur in the slip planes as a consequence of the localization of plastic slip. The propagation of micro-cracks follows crystallographic directions \cite{Rovinelli2018}, 
whereby it is an appropriate choice to couple damage initiation with crystal plasticity on the grain scale, cf. \cite{Flouriot2003}, \cite{Aslan2011}.

While fracture of grains has in the past often been modeled by continuum damage combined with crystal plasticity models, e.g. \cite{Aslan2011}, \cite{Ekh2004}, phase-field methods have become a popular choice for modeling ductile fracture coupled to plasticity, see for example the overview of different ductile phase-field models by \citeauthor{Alessi2018} \cite{Alessi2018} and references therein. 
One of the most advantageous features of phase-field modeling is that it allows arbitrary crack paths without suffering from a mesh dependence. 
On the micro-scale, phase-field models have been coupled to various crystal plasticity models, compare e.g. \cite{HernandezPadilla2014}, \cite{DeLorenzis2016}, \cite{Maloth2023}.
We have recently suggested a thermodynamic modeling framework for finite strain gradient-enhanced crystal plasticity, which has so far been applied to single crystals \cite{Auth2024}.

The largest difference in terms of modeling polycrystals opposed to single crystals is the treatment of grain boundaries.
Models combining gradient-extended plasticity and cohesive zone fracture at the grain boundaries of polycrystals have been suggested in e.g. \cite{Dahlberg2013} and \cite{Yalcinkaya2019}. Such models can be used to study the interplay between local stress concentrations at grain boundaries, which are influenced by the choice of grain boundary condition, and intergranular fracture. \citeauthor{Spannraft2020GrainDecohesion} \cite{Spannraft2020GrainDecohesion} have used a similar approach, but included an explicit coupling between the damage of cohesive element and the grain boundary condition. Specifically, the grain boundary condition that accounts for slip direction mismatch and grain boundary normal approaches a micro-free condition as the damage of the cohesive element grows.

In this contribution, we base the model formulation on the framework presented in \cite{Auth2024} for damage initiation and growth of micro-cracks in grains. 
The extension to the modeling of polycrystals in this paper requires a formulation of grain boundary conditions that should give size effects in the stress-strain response, but also allow for crack propagation across grain boundaries to yield realistic crack patterns.  
Therefore, we propose a thermodynamically consistent formulation of grain boundary conditions with a coupling to the phase-field damage such that an increased damage results in a lower slip transfer resistance between the grains.
The crystal plasticity phase-field fracture model from \cite{Auth2024} is also extended by a volumetric-deviatoric split to account for the influence of stress triaxiality. Since both the gradient extended crystal plasticity and phase-field fracture model result in highly nonlinear problems and demand very fine spatial discretization, a grand challenge when simulating failure in polycrystals in 3D is the computational time.
Therefore, we also present a robust and efficient numerical solution scheme. Furthermore, we formulate the gradient extension of the crystal plasticity to only act on a scalar measure whereby the size-dependence behavior still is captured while possible influence of slip direction mismatch is disregarded.  This opens the way for feasible numerical studies of relations between microstructural characteristics and transgranular failure mechanisms observed in 3D polycrystals.

The paper is structured as follows: Section \ref{section:thermodynamics} first supplies a recapitulation of the thermodynamic modeling framework from \cite{Auth2024} and subsequently proposed prototype model is described in Section \ref{section:prototype}.
Section \ref{section:weak} introduces the field equations and 
the domains of the different fields with respect to the grain boundaries.
It introduces different thermodynamically motivated choices of micro-boundary conditions, as well as the coupling between micro-flexible boundary conditions and the damage field.
In Section \ref{section:numerical_experiments} numerical experiments on polycrystals are presented and analyzed. The experiments yield insights into: the behavior of different choices of micro-boundary conditions upon damage development, the possibility to employ the model for void coalescence predictions and its capability to produce three-dimensional crack fronts.
Concluding remarks, as well as an outlook to possible future applications of the modeling framework are given in Section \ref{section:conclusions}.
\section{Thermodynamic framework}
\label{section:thermodynamics}
The model for gradient extended crystal plasticity and phase-field fracture of polycrystals is formulated within the thermodynamic framework presented in \cite{Auth2024}. This section presents this framework in brief for completeness.

The multiplicative decomposition of the deformation gradient $\bi F$ into an elastic part $\bi F_{\rm e}$ and a plastic part $\bi F_{\rm p}$ is given by 
\begin{equation}
    \bi F=\bi F_{\rm e} \cdot \bi F_{\rm p} \,.
\end{equation}
The free energy $\Psi$ is assumed to depend on the elastic Cauchy-Green deformation gradient $\bi C_{\rm e}=\bi F_{\rm e}^{\rm T} \cdot \bi F_{\rm e}$,
the set of isotropic hardening variables $\setka$, where $k_\alpha$ is the  isotropic hardening variable on the $\alpha$-th of $n_\alpha$ slip systems,
a damage (phase-field) variable $d$, the state variable $q$,   and the set of spatial gradients $\setnablaka$, as well as the spatial gradient $\bi \nabla_0 d$, according to
\begin{equation}
    \Psi=\Psi\left(\bi C_{\rm e}, \, q, \, \setka, \, \setnablaka, \, d, \, \boldsymbol{\nabla}_0 d \, \right) \,.
    \label{eq:generic_free_energy}
\end{equation}
The dissipation inequality when adopting quasi-static and isothermal conditions is given by
\begin{equation}
    \label{eq:dissipation_inequality_start}
    \int_{\V0} \boldsymbol P: \dot{\bi F} \, \dV0-\int_{\V0} \dot{\Psi} \, \dV0 \geq 0 \,,
\end{equation}
where $\bi P$ is the first Piola-Kirchhoff stress and $\V0$ represents the initial domain with outer boundaries $\A0$. Subsequently, we also introduce the domains $\Vi$, which refer to the grains of a polycrystal such that $\Vi$ represents the i-th grain in a polycrystal with $n_\mathrm{grains}$ grains in total. Their respective grain boundaries are referred to as $\Ai$.

Using the standard Coleman-Noll arguments \cite{Coleman1963} 
the elastic second Piola-Kirchhoff stress $\bi S_{\rm e}$ is obtained as (see e.g. \cite{Simo1988})
\begin{equation}
    \bi S_{\rm e}=2 \, \frac{\partial \Psi}{\partial \bi C_{\rm e}} \,.
    \label{eq:elastic_2nd_PK_stress}
\end{equation}
By introducing the Mandel stress ${\boldsymbol{M}}_\mathrm{e}=\bi C_{\rm e} \cdot \bi S_{\rm e}$ and the plastic velocity gradient ${\boldsymbol{L}}_\mathrm{p}=\dot{\bi F}_{\rm p} \cdot {\bi F}_{\rm p}^{-1}$, as well as using the divergence theorem, the dissipation inequality can be rewritten as 
\begin{equation}
    \mathcal{D}
    =
    \int_\V0
    \left(
        {\boldsymbol{M}}_\mathrm{e}  
        :
        {\boldsymbol{L}}_\mathrm{p}
        +
        Q \, \dot{q}
        +
        \sumalpha
            \kappa_\alpha \, \dot{k}_\alpha 
        +
        Y_\mathrm{d} \, \dot{d}
    \right)
    \, \dV0
    +
    \sumgrains
    \int_\Ai
        \sumalpha
           \kappa_\alpha^{\Gamma}  \, \dot{k}_\alpha
    \, \dA0
       +
    \int_\A0
       Y_\mathrm{d}^{\Gamma} \, \dot{d}
    \, \dA0
    \geq
    0
    \,.
    \label{eq:reduced_dissipation_inequality}
\end{equation}
Therein, the following dissipative quantities were introduced
\begin{align}
    \label{eq:kappa_def}
    \kappa_\alpha
    &=
    -\frac{\partial\Psi}{\partial k_\alpha} 
    +\boldsymbol{\nabla}_0 \cdot \frac{\partial\Psi}{\partial\boldsymbol{\nabla}_0 k_\alpha} \,,
    &\kappa_\alpha^{\Gamma}
    &=
    - \boldsymbol{N}
    \cdot
    \frac{\partial\Psi}{\partial\boldsymbol{\nabla}_0 k_\alpha}
    \,,
    \\ 
    Y_\mathrm{d}
    &=
    - \frac{\partial\Psi}{\partial d} 
    + \boldsymbol{\nabla}_0 \cdot \frac{\partial\Psi}{\partial\boldsymbol{\nabla}_0 d} \,,
    &Y_\mathrm{d}^{\Gamma}
    &=
    - \boldsymbol{N}
    \cdot
    \frac{\partial\Psi}{\partial\boldsymbol{\nabla}_0 d} \,,
    \\
    Q &= -\frac{\partial \Psi}{\partial q} \,,
\end{align}
and $\boldsymbol{N}$ is the unit normal to the boundaries.

\section{Prototype crystal plasticity model}
\label{section:prototype}
In this section, we present a prototype gradient-extended crystal plasticity phase-field fracture model.
We will assume the following form of the free energy, where a split of the elastic free energy contribution into a tensile part $\hat{\Psi}_\mathrm{e}^{+}$ and a compressive part $\hat{\Psi}_\mathrm{e}^{-}$ is adopted,
\begin{equation}
    \label{eq:free_energy}
    \Psi
    =
    g_\mathrm{e} \left(d, \, \epsilon^\mathrm{p}\right) \, \hat{\Psi}_\mathrm{e}^{+}\left(\boldsymbol{C}_\mathrm{e} \right)
    +
    \hat{\Psi}_\mathrm{e}^{-}\left(\boldsymbol{C}_\mathrm{e} \right)
    +
    \hat{\Psi}_\mathrm{p}\left(\setka, \, \setnablaka \right)
    +
    {\Psi}_\mathrm{d}\left(d, \boldsymbol{\nabla}_0 d\right)
    \,.
\end{equation}
$\hat{\Psi}_\mathrm{p}$ and $\Psi_\mathrm{d}$ represent the free energy contributions associated with gradient crystal plasticity and phase-field fracture, respectively. Following, the three contributions to the free energy are defined.
\subsection{Hyper-elasticity}
The effective (undamaged) elastic part of the free energy $\hat{\Psi}_\mathrm{e}$ is assumed to be of Saint-Venant type and a volumetric-deviatoric energy split is employed, whereby the tensile and compressive parts of the elastic free energy are given by
\begin{equation}
    \hat{\Psi}_\mathrm{e}^{+}\left(\boldsymbol{E}_\mathrm{e} \right)
    =
    \frac{1}{2} \, K \left\langle
        \mathrm{tr}\left(\boldsymbol{E}_\mathrm{e}\right)
    \right\rangle_{+}^2
    +
    G \, \boldsymbol{E}_\mathrm{e}^\mathrm{dev} : \boldsymbol{E}_\mathrm{e}^\mathrm{dev}
    \quad \text{and} \quad
    \hat{\Psi}_\mathrm{e}^{-}\left(\boldsymbol{E}_\mathrm{e} \right)
    =
    \frac{1}{2} \, K \left\langle
        \mathrm{tr}\left(\boldsymbol{E}_\mathrm{e}\right)
    \right\rangle_{-}^2
    \,,
\end{equation}
where $\boldsymbol{E}_\mathrm{e}=\left(\boldsymbol{C}_\mathrm{e} - \boldsymbol{I}\right)/2$ is the Green-Lagrange strain tensor and $\left\langle \bullet \right\rangle_\pm = \left(\bullet \pm \left|\bullet\right|\right)/2$ denotes the positive and negative parts of $\bullet$.
$K$ and $G$ represent the bulk and shear modulus, respectively.
The second Piola-Kirchhoff stress $\boldsymbol{S}_\mathrm{e}$ is then derived from Equation \ref{eq:elastic_2nd_PK_stress} as
\begin{equation}
    \boldsymbol{S}_\mathrm{e}
    =
    g_\mathrm{e}\left(d, \epsilon^\mathrm{p}\right)
    \, 2 \,
    \frac{\partial \hat{\Psi}_\mathrm{e}^{+}}{\partial \boldsymbol{C}_\mathrm{e}}
    + 
    2 \, \frac{\partial \hat{\Psi}_\mathrm{e}^{-}}{\partial \boldsymbol{C}_\mathrm{e}}
\end{equation}
Ductile failure is modeled by employing the elastic degradation function $g_\mathrm{e}$ presented by \citeauthor{Ambati2015} \cite{Ambati2015}
\begin{equation}
\label{eq:degradation_function}
    g_\mathrm{e}\left(d, \, \epsilon^\mathrm{p}\right)
    = 
    \left(1 - d \right)^{2\left(\epsilon^\mathrm{p} / \epsilon^\mathrm{p}_\mathrm{crit}\right)^n}
    \,,
\end{equation}
where $\epsilon^\mathrm{p}$ is the accumulated plastic strain, which is represented by $q$ in the thermodynamic modeling framework above.
Notice that this formulation requires the development of plasticity as well as phase-field damage in order to degrade the material. This is expected behavior for ductile failure in metals and also leads to the recovery of a true elastic zone independently of the choice of surface energy functional.
The critical accumulated plastic strain $\epsilon^\mathrm{p}_\mathrm{crit}$ and the exponent $n$ determine how the development of plastic strain impacts the material degradation.
As shown by \cite{Ambati2015}, this degradation function formulation is thermodynamically consistent, i.e. $Q \ \dot{q} \geq 0$. 

\subsection{Gradient-enhanced crystal plasticity}
The effective plastic free energy $\hat{\Psi}_\mathrm{p}$ is chosen as
\begin{equation}
\label{eq:plastic_free_energy}
    \hat{\Psi}_{\rm p}
    =
    \frac{1}{2} \, \sumalpha H_\alpha \, k_\alpha^2
    +
    \frac{l_g^2}{2} \, 
    \sumalpha
    \sumbeta
        \boldsymbol{\nabla}_0 k_\alpha \cdot \boldsymbol{H}_{\alpha\beta} \cdot \boldsymbol{\nabla}_0 k_\beta
    \,, 
\end{equation}
where the following material parameters were introduced: 
the isotropic hardening modulus $H_\alpha$ and the length scale for gradient-enhanced hardening $l_\mathrm{g}$.
Additionally, $\boldsymbol{H}_{\alpha\beta}$ determines gradient-enhanced cross-hardening between slip systems $\alpha$ and $\beta$. In general, the model allows to make distinct choices of $\boldsymbol{H}_{\alpha\beta}$ for gradient-enhanced self-hardening ($\alpha=\beta$) and gradient-enhanced cross-hardening ($\alpha \neq \beta$). However, in this work we choose $\boldsymbol{H}_{\alpha\beta} = H^\mathrm{g} \, \boldsymbol{I}$ to be constant for all slip systems $\alpha,\,\beta$. 
$H^\mathrm{g}$ can then be seen as gradient-enhanced hardening modulus.
We adopt the standard assumption that the slip direction $\bar{\bi s}_\alpha$ and normal vector to the slip plane $\bar{\bi m}_\alpha $ on the intermediate  configuration are fixed (and equal to their corresponding vectors on the undeformed configuration).
The yield functions for the slip systems $\Phi_\alpha$ are defined in terms of the effective Schmid stresses $\hat{\tau}_\alpha$
\begin{equation}
   \Phi_\alpha
   =
   \left| \hat{\tau}_\alpha \right|
   -
   \left(
       \tau_{\rm y} + \kappa_\alpha
    \right) \,,
\end{equation}
with $\hat{\tau}_\alpha={\tau}_\alpha\,/\,g_\mathrm{e}\left(d, \, \epsilon^\mathrm{p}\right)$, wherein $\tau_\alpha$ is the standard crystal plasticity Schmid stress $\tau_\alpha=\bi{M}_{\rm e}^\mathrm{dev} : \left( \bar{\bi s}_\alpha \otimes \bar{\bi m}_\alpha \right)$.  $\tau_{\rm y}$ is the initial yield stress, which is the same for all slip systems here.
We chose an associative evolution equation for the plastic velocity gradient  
\begin{equation}
\label{eq:evolution_Fp}
    \bar{\boldsymbol{L}}_{\rm p}
    =
    \dot{\bi F}_{\rm p} \cdot \bi F_{\rm p}^{-1}
    =
    \sumalpha
        \dot{\lambda}_\alpha \, \frac{\partial \Phi_\alpha}{\partial {{\bi M}}_{\rm e}}
    =
    \sumalpha
        \frac{\dot{\lambda}_\alpha}{g_\mathrm{e}\left(d, \, \epsilon^\mathrm{p}\right)} \, \left(
            \bar{\boldsymbol{s}}_\alpha \otimes \bar{\boldsymbol{m}}_\alpha \,
        \right)
        \, \mathrm{sign}\left(\hat{\tau}_\alpha\right)
\end{equation}
and apply a viscoplastic regularization for the multiplier $\lambda_\alpha$
\begin{equation}
    \label{eq:lambda_alpha}
    \dot{\lambda}_\alpha = \frac{1}{ t^\ast}\, \left< \frac{\Phi_\alpha}{\sigma_\text{d}} \right> ^ {m} \,,
\end{equation}
where $t^\ast$, $m$ and $\sigma_\text{d}$ control the viscosity of the model and $<\bullet>$ denotes Macauley brackets.
The accumulated plastic strain $\epsilon^{\rm p}$ is defined in terms of $\dot{\lambda}_\alpha$ 
\begin{equation}
\label{eq:evolution_ep}
    \epsilon_{\rm p}
    =
    \int_0^t 
    \sqrt{
        \sumalpha \dot{\lambda}_\alpha^2
    }
    \, {\rm d}t \,.
\end{equation}
The isotropic and gradient-extended hardening stresses $\kappa_\alpha$ and their boundary tractions $\kappa_\alpha^\Gamma$, respectively, are derived from Equation (\ref{eq:kappa_def}) as
\begin{align}
    \label{eq:prototype_kappa}
    \kappa_\alpha
    &=
    -H_\alpha \, k_\alpha
    +
    H^\mathrm{g} \, l_\mathrm{g}^2 \,
    \sumbeta
        \left(
            \boldsymbol{\nabla}_0 \cdot \boldsymbol{\nabla}_0 k_\beta
        \right) \,,
    \\
    \label{eq:prototype_kappagamma}
    \kappa_\alpha^\Gamma
    &=
    -H^\mathrm{g} \, l_\mathrm{g}^2 \,
    \boldsymbol{N}
    \cdot
    \sumbeta
         \boldsymbol{\nabla}_0 k_\beta \,.
\end{align}
An associative evolution is assumed for the hardening strains $k_\alpha$
\begin{equation}
\label{eq:evolution_k}
    \dot{k}_\alpha
    =
    \dot{\lambda}_\alpha \, \frac{\partial \Phi_\alpha}{\partial \kappa_\alpha}=-\dot{\lambda}_\alpha \,.
\end{equation}
These assumptions for the hardening can be extended, see e.g. \cite{Bargmann2010} to account for more complex models such as kinematic hardening, cross-hardening and nonlinear hardening.

\subsection{Phase-field fracture model}
An AT2 surface energy functional $\Gamma_\mathrm{d}$ (cf. \cite{Ambrosio1990}) is adopted for the phase-field model. The corresponding free energy contribution is given by
\begin{equation}
 {\Psi}_\mathrm{d}\left(d, \boldsymbol{\nabla}_0 d\right)=
 \mathcal{G}^\mathrm{d}_0 \, \Gamma_\mathrm{d}\left(d, \boldsymbol{\nabla}_0 d \right) \,
 \quad
 \mbox{with} \;
   \Gamma_\mathrm{d}  =
    \frac{1}{2 \, \ell_0}
    \left(
        d^2 + \ell_0^2 \left| \boldsymbol{\nabla}_0 d \right|^2
    \right)
\end{equation}
where $\mathcal{G}^\mathrm{d}_0$ represents fracture toughness and $\ell_0$ is the length-scale parameter controlling the width of the diffuse crack model.
\subsubsection{Damage irreversibility}
\label{section:irreversibility}
Since phase-field models do not inherently possess damage irreversibility, it is necessary to incorporate a model extension in order to avoid damage healing upon unloading. One of the most common choices is the history variable approach \cite{Miehe2010_historyvariable}. It has, however, been shown to be variationally inconsistent. 
Within this work we instead employ a micromorphic approach \cite{Forest2009}, which has recently been shown to allow for a variationally consistent framework that locally enforces damage irreversibilty \cite{Bharali2023}. 
In \cite{Auth2024}, we showed that this approach works well together with a gradient-enhanced crystal plasticity model.
In the micromorphic irreversibility formulation, an additional local variable $\varphi$ is introduced, which represents the phase-field damage in the quadrature points. The local phase-field damage $\varphi$ and the global phase-field damage $d$ are then connected by a penalty term in the free energy, that is controlled by a penalty parameter $\alpha$.
All occurrences of $d$ in the free energy, apart from that in the new penalty term and the gradient term $\boldsymbol{\nabla}_0 d$, are then replaced by the new local damage variable $\varphi$.
\begin{equation}
    \label{eq:free_energy_micromorphic}
    \Psi
    =
    g_\mathrm{e} \left(\varphi, \, \epsilon^\mathrm{p}\right) \, \hat{\Psi}_\mathrm{e}^{+}\left(\boldsymbol{C}_\mathrm{e} \right)
    +
    \hat{\Psi}_\mathrm{e}^{-}\left(\boldsymbol{C}_\mathrm{e} \right)
    +
    \hat{\Psi}_\mathrm{p}\left(\setka,\, \setnablaka\right)
    +
    {\Psi}_\mathrm{d}\left(\varphi, \boldsymbol{\nabla}_0 d\right)
    +
    \frac{\alpha}{2} \left( \varphi-d\right)^2
\end{equation}
The dissipation inequality (\ref{eq:reduced_dissipation_inequality}) is thereby modified to
\begin{equation}
    \mathcal{D} =
        \int_\V0
        \left(
                {\boldsymbol{M}}_\mathrm{e}  
                :
                {\boldsymbol{L}}_\mathrm{p}
                + Q \, \dot{q}
   + 
        \sumalpha
           \kappa_\alpha \, \dot{k}_\alpha
        +Y_\varphi\, \dot{\varphi}+
           Y_\mathrm{d}\, \dot{d}
        \right)
        \, \dV0
    +
    \sumgrains
        \int_\Ai
            \sumalpha
                \kappa_\alpha^{\Gamma}  \, \dot{k}_\alpha
        \, \dA0
       +
       \int_\A0
       Y^{\Gamma}
            \, \dot{d}
        \, \dA0
    \geq
    0 \,,
    \label{eq:diss_inequality_micromorphic}
\end{equation}
where $Y_\varphi = - \partial \Psi / \partial \varphi$. 
Assuming that $Y_\mathrm{d}=0$ ($Y_\mathrm{d}$ energetic) yields the global phase-field equation
\begin{equation}
    Y_\mathrm{d}
    =
    \alpha \, (\varphi-d)
    + \mathcal{G}_0^\mathrm{d}\,\ell_0 \,\boldsymbol{\nabla}_0 \cdot \boldsymbol{\nabla}_0 d
    = 0 \,
\label{eq:pf-micromorphic}
\end{equation}
The penalty term introduces an additional regularization to the model. In \cite{Miehe2017} it was mainly used for numerical robustness, but here it is used to ensure damage irreversibility.

In order to exploit the micromorphic formulation for damage irreversibility, the evolution of the local phase-field $\varphi$ is derived from the corresponding term of the dissipation inequality: $Y_\varphi \, \dot{\varphi} \geq 0$. For the proposed prototype model, $Y_\varphi=-\partial \Psi / \partial \varphi$ is obtained as 
\begin{equation}
\label{eq:local_residual_pf}
    Y_\varphi
    =
    - \frac{\partial g_\mathrm{e}}{\partial \varphi}
    \hat{\Psi}_\mathrm{e}^{+}
    -
    \frac{\mathcal{G}_0^\mathrm{d}}{\ell_0} \,\varphi 
    -
    \alpha \, \left(\varphi-d\right) \,.
\end{equation}
In the first step a trial damage value $\varphi^\mathrm{trial}$ is computed such that $Y_\varphi\left(\varphi^\mathrm{trial}\right)=0$. The trial value is then compared to the previous damage value in order to decide if it is accepted. Analogously to plasticity algorithms, we can formulate Karush-Kuhn-Tucker conditions for this process
\begin{equation}
    \label{eq:varphi_KKT}
    \dot{\varphi} \geq 0 \,, \quad 
	    \dot{\varphi} \, f_\mathrm{\varphi}=0 \,, \quad  
	    f_\mathrm{\varphi} = \varphi^\mathrm{trial}-\varphi \leq 0 \,.
\end{equation}
Thereby, thermodynamic consistency is fulfilled. Notice that the local phase-field $\varphi$ also is a history variable in this formulation. The authors would also like to point out that the micromorphic phase-field irreversibility is a self-contained feature of the model formulation and can easily be exchanged for a more classical history variable irreversibility approach. In this case the phase-field equation is still obtained from assuming that $Y_\mathrm{d}=0$, but $Y_\mathrm{d}$ is derived from Equation (\ref{eq:free_energy}) instead of Equation (\ref{eq:free_energy_micromorphic}).

\section{Weak form of balance equations}
\label{section:weak}
We express the weak form of the equilibrium equation in terms of the first Piola-Kirchhoff stress $\boldsymbol{P}$ and neglect inertial forces as well as body forces
\begin{equation}
    \delta\mathcal{W}^\mathrm{u}
    =
   \int_\V0
        \boldsymbol{P} : \left( \delta \boldsymbol{u} \otimes 
        \boldsymbol{\nabla}_0 \right)
    \, \dV0
    -\int_\A0
        \boldsymbol{t}_0^\ast
        \cdot \delta \boldsymbol{u}
    \, \dA0
    \,,
\label{eq:weak_equil}
\end{equation}
where  $\boldsymbol{t}_0^\ast$ is a prescribed traction on the boundary $\partial V_0$. 
Since the expression for the gradient-extended hardening stress $\kappa_\alpha$, Equation (\ref{eq:kappa_def}), includes spatial gradients, it also is a field equation. However, we adopt a dual mixed procedure \cite{Ekh2007}, whereby a gradient hardening field $\boldsymbol{g}$ is defined as
\begin{equation}
    \boldsymbol{g}
    =
    \sumalpha
        \boldsymbol{\nabla}_0 k_\alpha \,.
    \label{eq:strong_g}
\end{equation}
The weak form of Equation (\ref{eq:strong_g}) is then (instead of Equation (\ref{eq:kappa_def})) introduced as a field equation and by using the divergence theorem we obtain the weak form for gradient (crystal) plasticity as
\begin{equation}
    \delta\mathcal{W}^g_i
    =
    \int_{\Vi} 
         \boldsymbol{g} \cdot \delta  \boldsymbol{g}
    \, \dV0
    +
    \int_{\Vi}
        \sumalpha k_\alpha \, \boldsymbol{\nabla}_0 \cdot \delta \boldsymbol{g}
    \, \dV0
    -
    \int_\Ai
        \sumalpha k_\alpha \, \boldsymbol{N} \cdot \delta \boldsymbol{g}
    \, \dA0
    \quad \forall \: i = \{1,\,...\,,\,n_\mathrm{grains}\}
    \,.
\label{eq:weak_galpha}
\end{equation}
The gradient-extended hardening stress $\kappa_\alpha$ (Equation \ref{eq:prototype_kappa}) then becomes a local equation. 
\begin{equation}
    \kappa_\alpha
    =
    -H_\alpha \, k_\alpha
    +
    H^g \, l_g^2 \,
    \boldsymbol{\nabla}_0 \cdot \boldsymbol{g}
    \,,
\end{equation}
This procedure has been shown to be numerically robust (\cite{Ekh2007}, \cite{Carlsson2017}). Additionally, it allows to easily exchange the underlying local crystal plasticity model, since such changes only affect local equations..
Finally, the weak form of the global phase-field equation (\ref{eq:pf-micromorphic}) is given by
\begin{equation}
\delta\mathcal{W}^\mathrm{d}=
\int_\V0
        \alpha \left(\varphi - d\right) \, \delta d
    \, \dV0
    +
     \int_\A0
        {\,\mathcal{G}_0^\mathrm{d}\,\ell_0} \, \bi{N} \cdot 
        \boldsymbol{\nabla}_0 d \,  \delta d
    \, \dA0
    -
    \int_\V0
        {\,\mathcal{G}_0^\mathrm{d}\,\ell_0} \, \boldsymbol{\nabla}_0 d \cdot \boldsymbol{\nabla}_0 \delta d
    \, \dV0
\label{eq:weak_pf}
\end{equation}
where the standard boundary condition $\bi{N} \cdot 
\boldsymbol{\nabla}_0 d =0$ will be assumed on the outer boundary of the grain structure, such that the boundary integral term disappears.
The micromorphic penalty parameter should be chosen in relation to the effective fracture energy $\mathcal{G}^\mathrm{d}_0/\ell_0$, e.g. by employing a dimensionless scalar $\beta$ and chosing $\alpha = \beta \, \mathcal{G}^\mathrm{d}_0 / \ell_0$.

Equations (\ref{eq:weak_equil}), (\ref{eq:weak_galpha}) and (\ref{eq:weak_pf}) are to be solved for polycrystals, which requires the definition of how the respective fields should behave at grain boundaries. Since this work is concerned with transgranular fracture, the displacement field $\boldsymbol{u}$ is assumed to be continuous across grain boundaries. The gradient field $\boldsymbol{g}$ on the other hand is a measure for the dislocation distribution and since grain boundaries restrict the dislocation movement in polycrystals, the gradient field $\boldsymbol{g}$ is in general discontinuous. While dislocation transport across grain boundaries can occur in reality, we here make the simplification that the gradient field $\boldsymbol{g}$ is independent in each grain.
The damage field $d$ is assumed to be continuous across grain boundaries.
\subsection{Micro-flexible boundary conditions}
\label{sec:microflexible_bcs}
Requirements on the boundary conditions for the strain gradient field $\boldsymbol{g}$ follow from thermodynamic consistency. For the chosen prototype model, the boundary term relating to gradient (crystal) plasticity in the reduced dissipation inequality, Equation (\ref{eq:reduced_dissipation_inequality}), is given by
\begin{equation}
    \label{eq:dissipation_inequality_galpha_boundary}
    \sumgrains
    \int_\Ai
        \sumalpha
            \left(
                \kappa^\Gamma_\alpha \, \dot{k}_\alpha 
            \right)
    \, \dA0
    =
    \sumgrains
    \int_\Ai
        \left(
            -H^\mathrm{g} \, l_\mathrm{g}^2 \,
            \boldsymbol{N}
            \cdot
            \boldsymbol{g}
        \right)
        \sumalpha
            \dot{k}_\alpha
    \, \dA0
    \geq 0 \,.
\end{equation}
This is a stronger requirement than postulated by Equation (\ref{eq:reduced_dissipation_inequality}), however we will prove that the proposed boundary conditions can satisfy this requirement and thereby also satisfy Equation (\ref{eq:reduced_dissipation_inequality}).
There are two trivial ways to fulfill the inequality in Equation (\ref{eq:dissipation_inequality_galpha_boundary}):
\begin{enumerate}
    \item Micro-free boundary conditions: $\kappa_\alpha^\Gamma=0$, corresponding to Dirichlet-type constraints by prescribing $\boldsymbol{N} \cdot \boldsymbol{g} = 0$ in the dual-mixed formulation.
    \item Micro-hard boundary conditions: $\sumalpha k_\alpha = 0$ and thereby $\dot{k}_\alpha = 0$. This corresponds to removing the boundary integral in the weak form for gradient plasticity, Equation (\ref{eq:weak_galpha}).
\end{enumerate}
Additionally to these two classical types of boundary conditions, we propose micro-flexible boundary conditilns inspired by the formulation in \cite{Ekh2011}.
We prescribe the sum of hardening strains on the boundary such that
\begin{equation}
    \sumalpha k_\alpha = C_\Gamma \, \kappa^\Gamma \,,
\end{equation}
where a micro-flexibility parameter $C_\Gamma$ has been introduced and 
\begin{equation}
    \kappa^\Gamma = \kappa^\Gamma_\alpha = -H^\mathrm{g} \, l_\mathrm{g}^2 \, \boldsymbol{N} \cdot \boldsymbol{g} \,,   
\end{equation}
since $\kappa^\Gamma_\alpha$ for the proposed model is the same for all slip systems.
The boundary dissipation in Equation (\ref{eq:dissipation_inequality_galpha_boundary}) thereby becomes
\begin{equation}
    \sumgrains
    \int_\Ai
        \frac{1}{C_\Gamma} \, 
        \sumalpha
            k_\alpha
        \,
        \sumalpha
            \dot{k}_\alpha
    \, \dA0
    \geq 0 \,.
\end{equation}
Therein, all signs are known: From the evolution equations, Equation (\ref{eq:evolution_k}), we know that $\dot{k}_\alpha \leq 0$ and with an initial value of $k_\alpha^\mathrm{initial} = 0$, it follows that $k_\alpha \leq 0$. Thus, by choosing $C_\Gamma > 0$ thermodynamic consistency of the boundary term is guaranteed.
Upon applying micro-flexible boundary conditions the weak form for gradient plasticity, Equation (\ref{eq:weak_galpha}) becomes
\begin{equation}
    \label{eq:weak_galpha_microflexible}
    \delta\mathcal{W}^g_i
    =
    \int_{\Vi} 
         \boldsymbol{g} \cdot \delta  \boldsymbol{g}
    \, \dV0
    +
    \int_{\Vi}
        \sumalpha k_\alpha \, \boldsymbol{\nabla}_0 \cdot \delta \boldsymbol{g}
    \, \dV0
    +
    \int_{\Ai} 
        \left(
            C_\Gamma \,
            H^\mathrm{g} \, l_\mathrm{g}^2 \, \boldsymbol{N} \cdot \boldsymbol{g}
        \right)
         \,
        \left( 
             \boldsymbol{N} \cdot \delta \boldsymbol{g}
         \right)
    \, \dA0
    \quad \forall \: i = \{1,\,...\,,\,n_\mathrm{grains}\}
    \,.
\end{equation}
The boundary term therein can be implemented as a standard boundary integral (over all grain boundaries) in the finite element framework, as it only includes constants, the normal vector to the boundary $\boldsymbol{N}$ and the field variable $\boldsymbol{g}$.
Note that for $C_\Gamma \rightarrow 0$ micro-hard boundary conditions are recovered. 
Micro-free boundary conditions are recovered when $C_\Gamma\rightarrow \infty$, since in this case $\kappa^\Gamma \rightarrow 0$ is required for satisfying the residual equation.
Grain boundaries present obstacles for the propagation of plastic slip, thereby causing a grain size dependent stress-strain response in polycrystals, also known as Hall-Petch the effect. The resistance to plastic slip across a grain boundary depends, e.g., on the misorientation angle between the slip systems of the grains meeting at the boundary and their accumulated plastic slip. See formulations for this in e.g. \cite{Gurtin2008}, \cite{Ekh2011},  and \cite{McBride2016}.
Varying resistance to plastic slip from grain boundary to grain boundary can be accounted for by adjusting the (initial) value of the micro-flexibility parameter $C_\Gamma$.
However, the dependence of the misorientation angles between single slip systems is for simplicity disregarded in this work in order to focus on transgranular fracture and the interaction with slip transmission resistance as well as to obtain feasible computational times.
In order to account for the effect of fracture on the resistance of grain boundaries to plastic slip, we propose to couple the micro-flexibility parameter $C_\Gamma$ with the  phase-field damage $d$ such that an increased $d$ lowers the resistance. In this way a boundary can evolve from being initially micro-hard to a micro-free boundary. For the presented prototype model, a linear relation between the micro-flexibility and the phase-field damage is chosen according to
\begin{equation}
    C_\Gamma\left(d\right) = C_{\Gamma,0} + C_\Gamma^\mathrm{d} \, d \,,
\end{equation}
where $C_{\Gamma,0}$ is the initial micro-flexibility of each grain boundary and $C_\Gamma^\mathrm{d}$ represents the additional amount of micro-flexibility gained by complete material degradation. The micro-flexibility parametrization should always be chosen relative to the gradient hardening parameters.
By this choice of $C_\Gamma$, the finite element problem, Equation (\ref{eq:weak_galpha_microflexible}), and the $\boldsymbol{g}$ -field become uncoupled between the grains (apart from an implicit coupling via $\boldsymbol{u}$ and $d$). 
\section{Numerical experiments}
\label{section:numerical_experiments}
The behavior of the proposed model is discussed with respect to different aspects of damage initiation in this section. First, examples of the base cases of micro-free and micro-hard boundary conditions are shown in order to discuss the qualities and flaws of both choices.
It is then demonstrated how micro-flexible boundary conditions coupled to the damage development combine the desirable features of those boundary conditions.
Finally, we show that the model is able to capture relevant phenomena in (ductile) fracture such as void coalescence and three-dimensional crack fronts.

\begin{figure}[htb]
    \centering
    \includegraphics{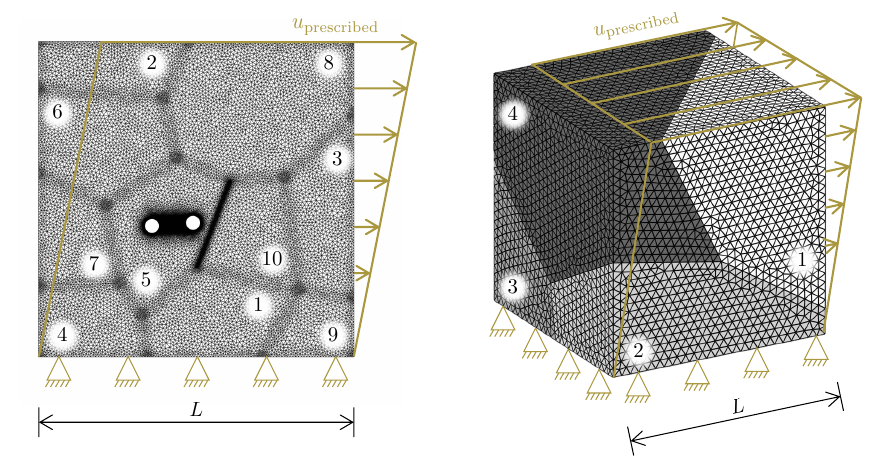}
    \caption{Two- (left) and three-dimensional (right) grain structures employed for the numerical experiments. The two-dimensional structure consists of 10 grains and is used with and without the voids in the middle grain. Without voids, the middle grain is meshed in the same manner as the other grains. The three-dimensional grain structure consists of four grains. All set-ups are loaded by simple shear along the horizontal axis.}
    \label{fig:meshes}
\end{figure}

Figure \ref{fig:meshes} shows the two different grain structures that are used for the numerical examples. The two-dimensional examples are performed on a structure consisting of 10 grains and employing plane strain conditions. The same grain structure is first investigated without the voids shown in Figure \ref{fig:meshes} and then with the voids. The three-dimensional examples are performed on a structure consisting of 4 grains. In all cases full FCC slip systems are used for all grains. The employed meshes contain mesh refinements along the grain boundaries and at grain boundary intersections. The two-dimensional meshes consist of $22\,616$ (without voids) and $45\,688$ (with voids) linear triangular elements, respectively. The three-dimensional mesh consists of $76\,506$ linear tetrahedral elements. Further details of the geometries and their meshing, are given in Appendix \ref{appendix:geometry_meshing}.
All structures are loaded in-plane by displacement-controlled pure shear loading.  
As introduced in Section \ref{section:weak}, two vector fields (displacement field $\boldsymbol{u}$ and gradient field $\boldsymbol{g}$) and a scalar field (global phase-field damage $d$) are needed on each node. While the displacement field $\boldsymbol{u}$ and the phase-field $d$ are continuous fields, the gradient field $\boldsymbol{g}$ is independent in each grain. This is achieved by duplicating the nodes on the grain boundaries and then applying linear constraints on the grain boundaries for the displacement field $\boldsymbol{u}$ and the phase-field $d$, which allows to use the same mesh representation for all fields.
Backward Euler time integration is applied to all time-dependent equations, see Appendix \ref{appendix:time_integration} for the respective derivations.
A staggered solver is used in order to solve the coupled equation system. The initial time step is $\Delta t_\mathrm{initial} = 0.1\,t^\ast$, the time step is halved whenever a load step does not converge. For an iterative solver to succeed on this type of phase-field problems, it is important that the staggered scheme extends to the local variables. Further, methods for avoiding unneeded iterations far away from the solution are employed. A detailed description of the solver algorithm and its respective parameters, is given in Appendix \ref{appendix:solver}.
This work focuses on model development rather than model parametrization, for details on material parameters, slip system orientations and loading, consult Appendix \ref{section:appendix_parameters}.

The implementation of the presented model is written in the \texttt{Julia} programming language \cite{Bezanson2017} and uses the \texttt{Ferrite.jl} finite element toolbox \cite{Ferrite.jl}. The material model implementation makes heavy use of the \texttt{Tensors.jl} tensor calculus toolbox \cite{Tensors.jl}, that supports automatic differentiation of tensor equations, and all figures are generated by the \texttt{Makie.jl} data visualization ecosystem \cite{Danisch2021}. 
\subsection{Impact of boundary conditions for gradient plasticity}
\begin{figure}[htb]
    \centering
    \includegraphics[height=10.3cm]{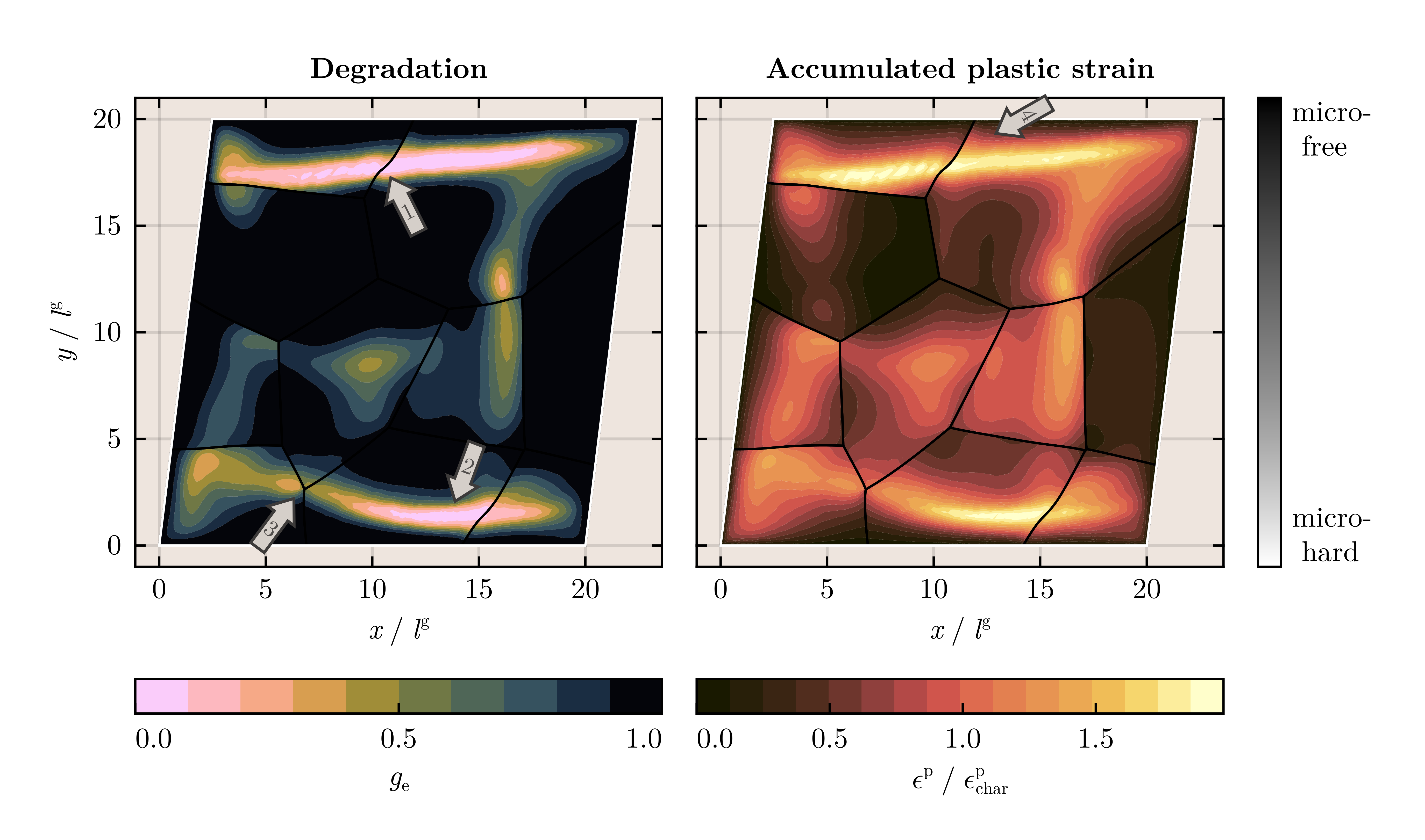}
    \caption{Degradation and plastic strain response for \textbf{micro-free boundary conditions} on the inner grain boundaries. Both fully developed cracks cross grain boundaries. Localization of plastic strain initially occurs in the grain marked with arrow 1, close to the grain boundary. The bottom crack (arrow 2) develops later in the simulation. Its crack growth direction is impacted by the stress concentration caused by the grain boundary intersection marked with arrow 3. The micro-free boundary conditions are reflected by the accumulated plastic strain contour lines that are perpendicular to the grain boundaries, especially in the highly plastified regions (e.g. arrow 4). Micro-hard boundary conditions are prescribed on the outer boundaries of the grain structure for numerical stability.}
    \label{fig:microfree}
\end{figure}
\begin{figure}[htb]
    \centering
    \includegraphics[height=10.3cm]{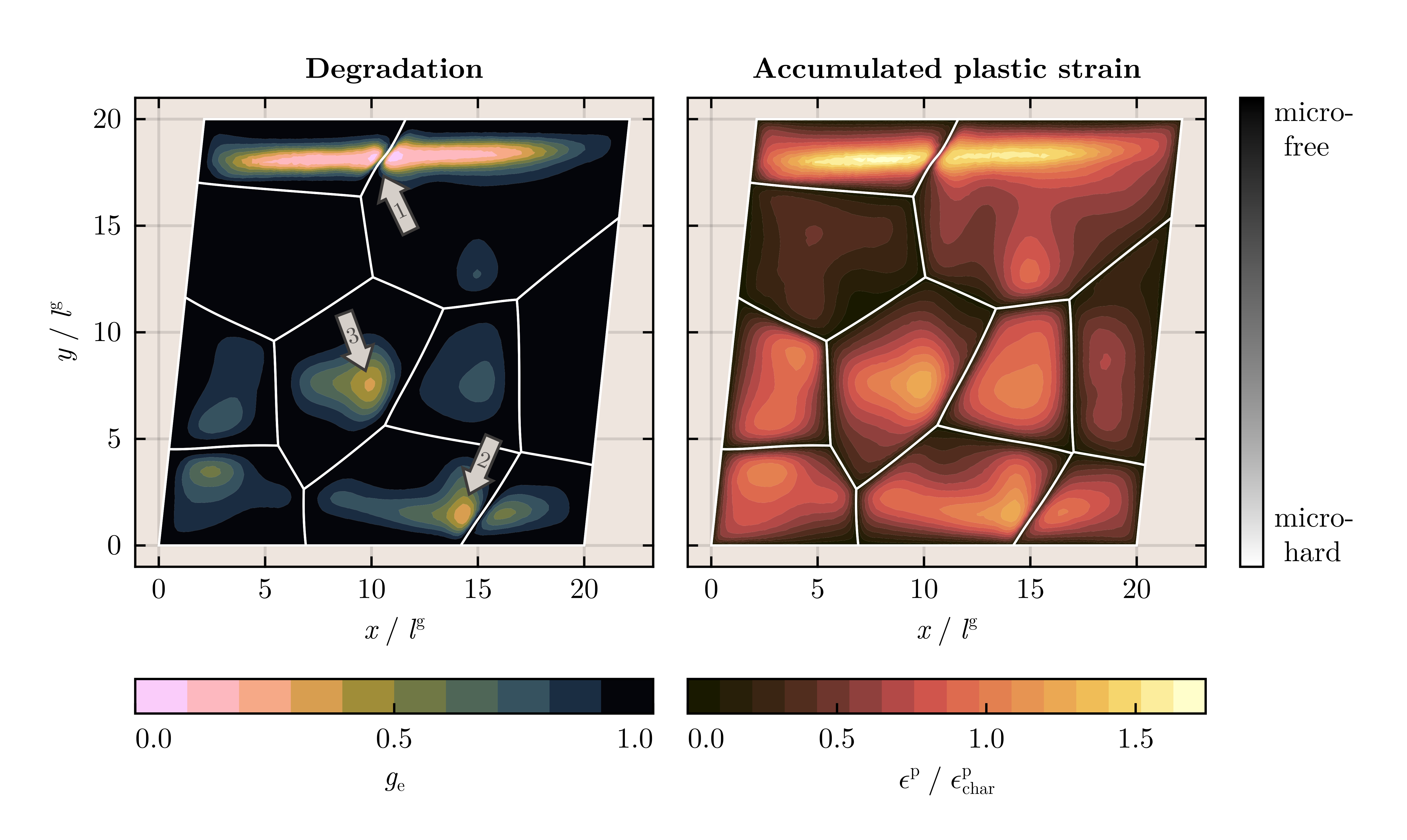}
    \caption{Degradation and plastic strain response for \textbf{micro-hard boundary conditions} on the inner grain boundaries. The crack developing within the upper grains is interrupted on the grain boundary marked by arrow 1. Localization of plastic strain initially occurs inside the grain right of that grain boundary. Plastic strain cannot develop on the grain boundaries. Micro-hard boundary conditions are prescribed on the outer boundaries of the grain structure.}
    \label{fig:microhard}
\end{figure}
\textbf{Micro-free and micro-hard boundary conditions} \newline
Figures \ref{fig:microfree} and \ref{fig:microhard} show the bounding cases of micro-free and micro-hard boundary conditions on the inner grain boundaries. Note that micro-hard boundary conditions are employed for the outer boundaries in all simulations in order to stabilize the global response for the Newton-Raphson solver. To the left of both figures, the degradation $g_\mathrm{e}$ is presented, displaying the cracks formed at the end of the simulation. To the right, the corresponding accumulated plastic strain $\epsilon^\mathrm{p}$ distribution, normalized by the characteristic accumulated plastic strain $\epsilon^\mathrm{p}_\mathrm{char}$ is shown. In both figures, the grain boundaries are colored in gray scale where dark means micro-free and light means micro-hard boundary conditions.

The accumulated plastic strain is an important input to the degradation function, as well as a representative quantity for the hardening variables $k_\alpha$.
Micro-free boundary conditions require that the component of the gradient hardening field $\boldsymbol{g}$ perpendicular to the  grain boundaries is zero. In Figure \ref{fig:microfree} this is reflected by contour lines of $\epsilon^\mathrm{p}$ which are perpendicular to the boundaries, especially on those grain boundaries where a large amount of plastic strain has developed (e.g. pointed out by arrow 4).
Notice however, that $\boldsymbol{g} = \sumalpha \boldsymbol{\nabla}_0 \, k_\alpha$ is only fulfilled in a weak sense. 
Micro-free boundary conditions allow the development of plastic strain on the grain boundaries. Initial strain localization is to a large extent influenced by the alignment of crystal orientations with the loading direction, as well as by geometrical features. Initial localization sites are often located at grain boundaries or grain boundary intersections. The upper crack for example initiates on the right side of the grain boundary marked by arrow 1. Later in the simulation a second crack develops, which initiates in the grain marked by arrow 2. The growth direction of the secondary crack is strongly influenced by the stress concentration occurring at the grain boundary intersection marked by arrow 3. 
Micro-hard boundary conditions on the other hand require that $\sumalpha k_\alpha = 0$ on the boundary, in a weak sense. This is reflected in Figure \ref{fig:microhard} by contour lines of $\epsilon^\mathrm{p}$ which are parallel to the grain boundaries. In particular, micro-hard boundary conditions disallow the development of plastic strain on the grain boundaries.
This is important in combination with the choice of degradation function, see Equation (\ref{eq:degradation_function}), since it also results in disallowing damage development on the grain boundaries. Consequently, the first crack in Figure \ref{fig:microhard} (arrow 1), is interrupted by the respective grain boundary. The micro-hard simulation shows two damage initiation sides in addition to the developed crack, which are marked with arrows 2 and 3 in Figure \ref{fig:microhard}. Comparing with the crack pattern arising from the micro-free simulation, it can be observed that both simulations develop a similar primary crack. The precise shape of the crack is impacted by the choice of micro-free / micro-hard boundary conditions. Additionally, the damage initiation site marked by arrow 2 in the micro-hard results corresponds to the damage initiation site of the secondary crack in the micro-free results (compare Figure \ref{fig:microfree}). Thus, the presented model predicts that the boundary conditions on gradient hardening have a significant impact on the crack growth across grain boundaries, but not on the damage initiation sites.

Both cases, micro-free and micro-hard boundary conditions, exhibit physically partly correct and partly incorrect behavior for the representation of ductile fracture in polycrystalline metals. Grain boundaries represent obstacles for dislocation movement. This is represented by micro-hard boundary conditions, causing an elastic response close to the grain boundaries (plastic deformation is caused by dislocation movement, thus preventing dislocation movement leads to an elastic response). This is also how the size-dependent stress–strain response in polycrystals is captured by gradient-plasticity models. 
Additionally, damage initiation commonly occurs inside the grains and not on the grain boundaries as implied by the micro-free results shown in Figure \ref{fig:microfree}.
However, the resistance of a given grain boundary to dislocation movement depends among other things on the misalignment between the slip system orientations on both sides of the boundary. Grain boundaries also do not present ultimate obstacles for crack growth, as implied by the micro-hard results shown in Figure \ref{fig:microhard}. 

\begin{figure}[htb]
    \centering
    \includegraphics[height=10.3cm ]{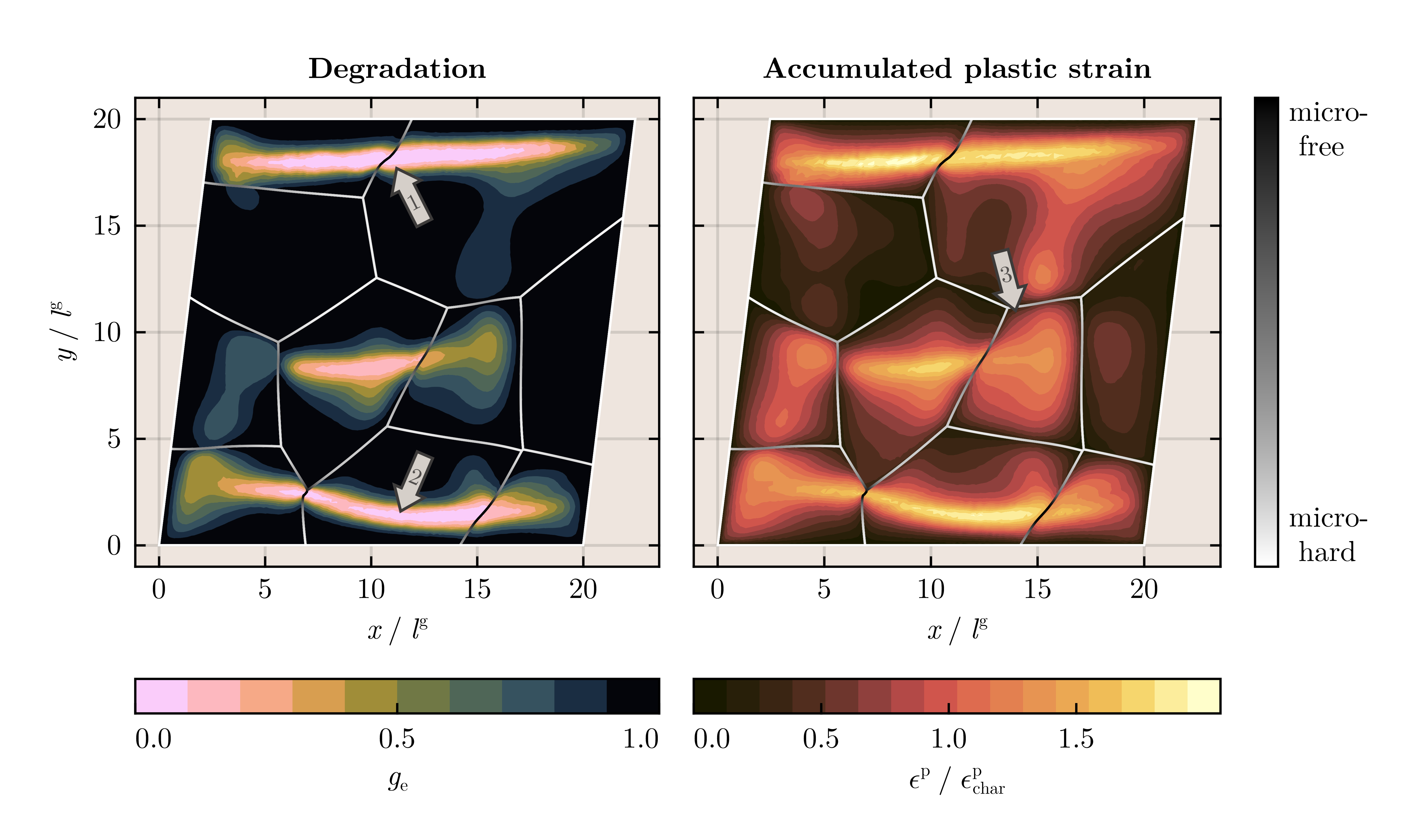}
    \caption{Degradation and plastic strain response for \textbf{micro-flexible boundary conditions} on the inner grain boundaries. 
    {The damage initiation side of the upper crack is indicated by arrow 1. The crack intiates on the right side of the grain boundary and then propagates to the left grain.}
    The micro-flexible boundary conditions initially behave nearly micro-hard and upon damage development locally transition to micro-free behavior. The state of the micro-boundary conditions at the end of the simulation is color-coded on a white to black color scale in the figure.
    Micro-hard boundary conditions are prescribed on the outer boundaries of the grain structure for numerical stability.}
    \label{fig:microflexible}
\end{figure}
\begin{figure}[htb]
    \centering
    \includegraphics{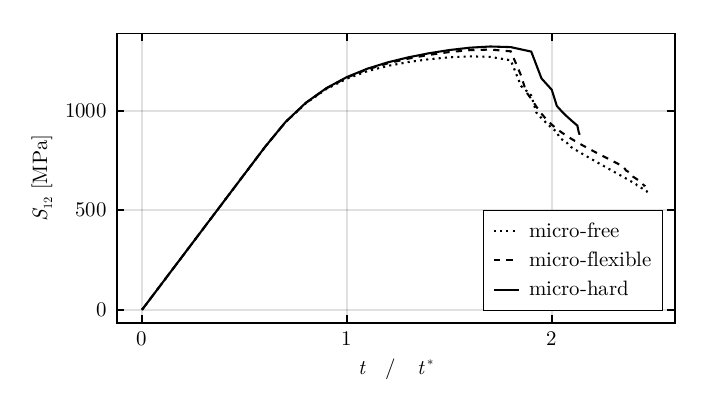}
    \caption{Volume averaged in-plane shear component of the second Piola-Kirchhoff stress tensor in the upper two grains, comparing different inner micro-boundary conditions. Micro-hard boundary conditions lead to a stiffer response than micro-free boundary conditions.
    As desired, micro-flexible boundary conditions initially behave similarly to micro-hard boundary conditions and then approach the behavior of micro-free boundary conditions during the softening period.
    } 
    \label{fig:S12}
\end{figure}
\textbf{Micro-flexible boundary conditions} \newline
In order to address the lack of accuracy of the micro-free and micro-hard boundary conditions in describing the full range from strain localization to damage initialization and finally fracture, we employ the micro-flexible boundary condition presented in Section \ref{sec:microflexible_bcs}. Under the assumption that the resistance of grain boundaries against dislocation movement decreases as damage develops, we suggest that the initial behavior of grain boundaries can be approximated by micro-hard boundary conditions while the behavior close to material failure should be approximated by micro-free boundary conditions. As described in Section \ref{sec:microflexible_bcs}, this is implemented by assigning a low initial micro-flexibility $C_{\Gamma,0}$ (low micro-flexibility recovers micro-hard behavior) that gradually increases as the global damage field $d$ increases (high micro-flexibility recovers micro-free behavior).
Figure \ref{fig:microflexible} displays the resulting fracture and plasticity patterns.
It can be seen that the obtained cracks each span several grains and pass the grain boundaries in a similar manner as for the micro-free case. By inspecting the white-to-black color scale on the grain boundaries, it can be noticed that the grain boundaries locally transition from micro-hard to micro-free behavior. It can be observed that most grain boundaries remain closer to the micro-hard state, while the boundary conditions locally change to micro-free at the locations where the material is most degraded. Arrow 3 marks a grain where the remaining micro-hard boundary conditions can be noticed particularly well from the contour lines of the accumulated plastic strain. Similarly to the micro-free case, the upper crack (arrow 1) develops first and the lower crack (arrow 2) develops later. The latter crack is also impacted by the same stress concentration as in the micro-free simulation.

Additionally, strain and damage localization are distinctly observable in Figure \ref{fig:microflexible}. 
The grain boundary marked by arrow 1 exhibits a clear s-shaped deformation, as well as a clear local transition to micro-free boundary conditions. This displays the need for a finite strain formulation, as well as the effect of the localized effect of micro-flexible boundary conditions.

Figure \ref{fig:S12} shows the volume-averaged shear component of the second Piola-Kirchhoff stress in the upper two grains (where the first crack develops in all simulations). It shows that the micro-hard boundary conditions lead to a stiffer response than the micro-free boundary conditions. Micro-flexible boundary conditions behave similarly to micro-hard boundary conditions during the hardening stage, but similarly to micro-free boundary conditions during the softening stage. Thus, they recover the physically meaningful behavior over the full plastification and fracture range.

\subsection{Void coalescence}
\begin{figure}[htb]
    \centering
    \includegraphics[height=10.3cm]{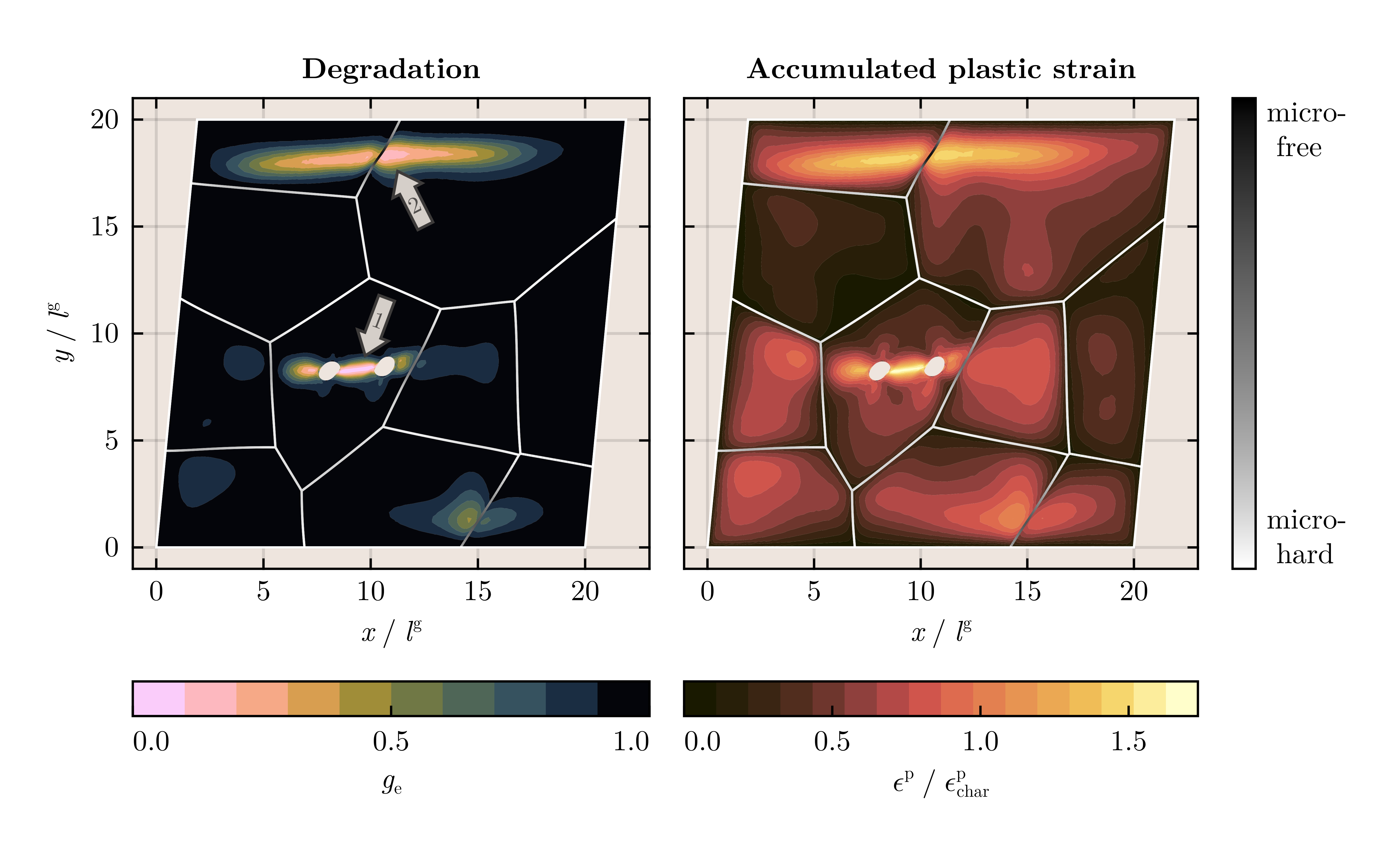}
    \caption{Degradation and plastic strain response for a case of void coalescence with \textbf{micro-flexible boundary conditions} on the inner grain boundaries. Pre-formed voids perturb the fracture pattern compared to Figure \ref{fig:microflexible} and lead to an additional crack, caused by void coalescence inside the middle grain (arrow 1). The multi-axial stress state around the voids activates slip systems that are otherwise not triggered by the global loading and thereby causes the failure of the originally undamaged grain.
    Micro-hard boundary conditions are prescribed on the outer boundaries of the grain structure for numerical stability.}
    \label{fig:microflexible_holes}
\end{figure}

Crack initiation often occurs due to void nucleation and growth.
Voids commonly nucleate on particles and inclusions placed within shear bands.
The following example shows damage growth from preexisting voids in a grain structure. The voids are here modeled as holes in the structure and micro-free boundary conditions are employed at the void boundaries.
The voids are positioned in the middle crack shown in Figure \ref{fig:microflexible}. This is the last crack that develops in the micro-flexible simulation. 
The resulting crack pattern is show in Figure \ref{fig:microflexible_holes}. 
The developed crack pattern is different from the one obtained in Figure \ref{fig:microflexible}. The cracks indicated by arrows 1 and 2 develop simultaneously, where the crack developing between the voids (arrow 1) is clearly more developed than the crack in the top grains (arrow 2).
Due to the crystal plasticity model the susceptibility of the grains to plasticity and thereby damage initiation highly depends on the crystal orientation of each grain relative to the loading direction. The multi-axial stress state caused by the voids leads to the activation of slip systems that were not triggered in the previous examples and thereby causes the failure marked by arrow 1.
This crack has a crack shielding effect, which prevents the development of the bottom crack compared to the original micro-flexible simulation.

\subsection{Three-dimensional crack fronts}
Crack fronts in polycrystals are strongly influenced by the grain structure, crystal orientations and the misalignment thereof, as well as voids and particles in the metal matrix. Hence, these crack fronts inherently propagate in three-dimensional space and it is not always possible to understand their behavior from two-dimensional simulations. The presented modeling framework is formulated in a dimension-agnostic manner and can thereby account for such three-dimensional effects. 
The final numerical example shows the response of a three-dimensional four-grained polycrystal to in-plane shear loading.
In order to allow for a coarser mesh, the length scale parameters have been adjusted to $l_\mathrm{g}=0.1333\,L$ and $\ell_0=0.08\,L$. The micro-flexiblity parameters, as well as the micromorphic penalty parameter are adjusted accordingly, see Table \ref{tab:material_parameters}.

\begin{figure}[htb]
    \centering
    \includegraphics[width=16cm]{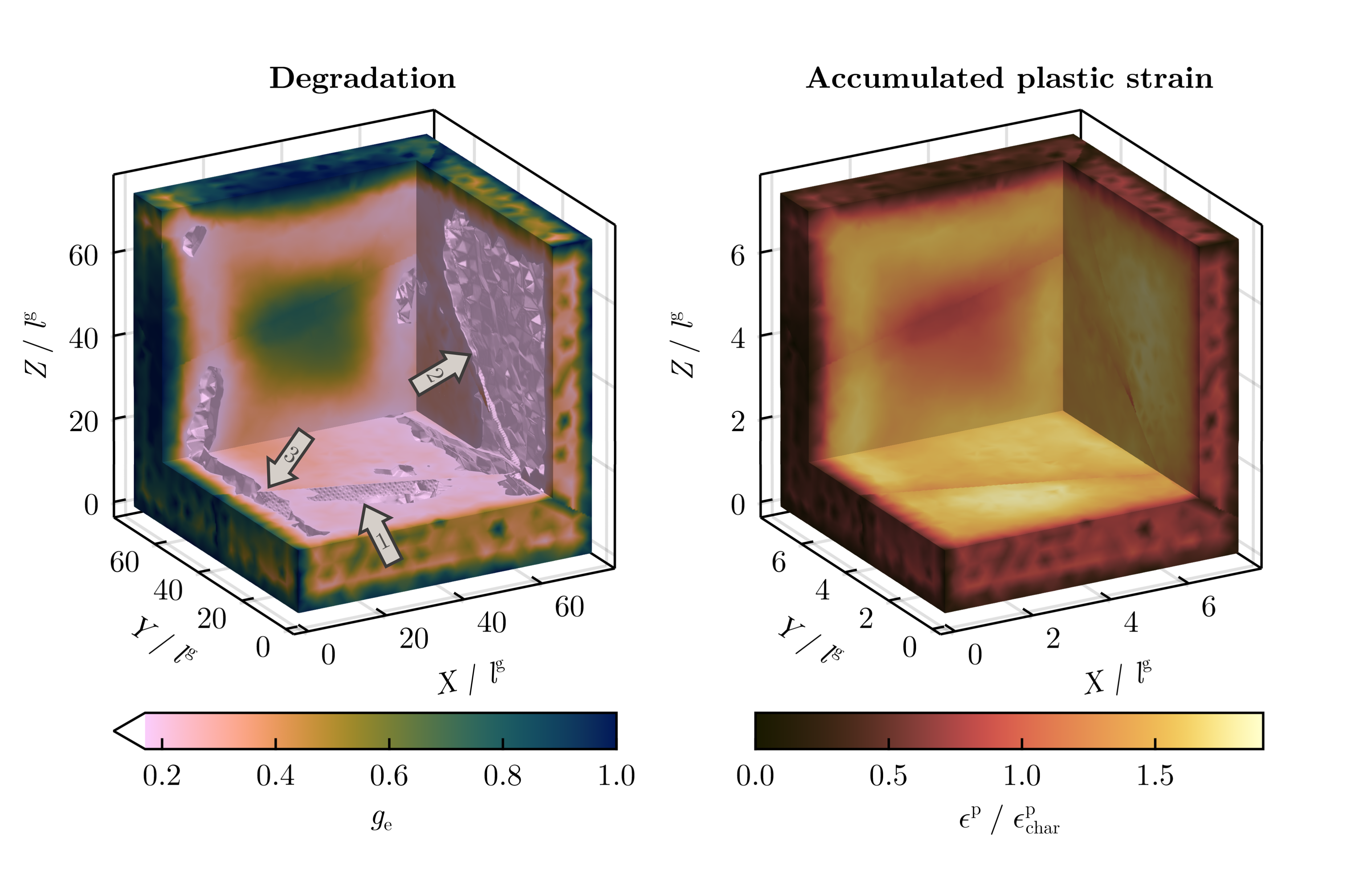}
    \caption{Degradation and plastic strain response for the four-grain three-dimensional grain structure. \textbf{Micro-flexible boundary conditions} are prescribed on the inner grain boundaries, while micro-hard boundary conditions are prescribed on the outer boundaries of the grain structure for numerical stability. 
    In the degradation plot all regions with degradation values $g_\mathrm{e} \leq 0.17$ are represented as transparent in order to better visualize the resulting crack fronts.
    The structure is oriented as presented in Figure \ref{fig:meshes} and thus shear is applied in the XZ-plane.
    Several crack fronts have developed. It can be observed that the crack fronts on the right and in the middle of the grain do not range along the entire Y-axis and could therefore not be captured by a plane strain set-up.
    }
    \label{fig:3d_results}
\end{figure}

Figure \ref{fig:3d_results} shows the resulting degradation and accumulated plastic strain fields. In order to better represent the crack fronts, all degradation values $g_\mathrm{e} \leq 0.17$ are represented as transparent. 
Simple shear is applied in the XZ-plane, with roller supports along the front and back faces in the Figure.
Notice that opposed to the previous two-dimensional results, the three-dimensional results are presented in the initial configuration.
The larger length scales for gradient plasticity and phase-field damage lead to wider shear bands and cracks. 
The final state of the simulation shows a complex crack front that has developed in the structure. The visualization in Figure \ref{fig:3d_results} shows a cut into part of the crack front. It can be observed that parts of the crack front do not range all the way along the Y-axis (marked by arrows 1 and 2), indicating that they could not have been resolved by a 2D plane strain simulation. This can in particular be observed at the right part of the crack front (arrow 2), which developed within grain 1 and does not cross the boundary to grain 4 (arrow 2 points on the grain boundary). 
The crack front developing on the left (arrow 3) on the other hand displays the impact of the micro-flexible boundary conditions. It reaches across the grain boundary between grain 2 and grain 3, the point where it crosses the grain boundary is marked by arrow 3.

\section{Concluding remarks}
\label{section:conclusions}
In this work, we have presented a modeling framework and a prototype model for ductile phase-field fracture across grain boundaries of polycrystals. The model fits a thermodynamic framework, that we have presented in previous work \cite{Auth2024} in the context of single crystal plasticity.

The prototype model presented in our previous work has been extended in order to account for properties of polycrystals. In particular, we have presented a damage-dependent micro-flexible boundary condition that allows for  modeling of grain boundaries that exhibit varying slip transmission resistance. The level of resistance is tunable by a micro-flexibility parameter, which we employ to locally vary boundary conditions from micro-hard to micro-free upon damage development. The presented boundary condition is represented by a standard boundary integral including field variables and material parameters. Hence, it is straightforward to implement it in a finite element setting.
Further, we have added a volumetric-deviatoric energy split to the prototype model in order to allow a physically reasonable difference in fracture behavior of tension and compression.
This is relevant for polycrystals, since the micro-stress state caused by the grain structure can easily yield local compression even under shear loading. It is also important for obtaining realistic behavior upon crack closure.

We have moreover employed a single gradient hardening field in the dual mixed formulation, instead of one gradient field per slip system as presented in \cite{Auth2024}. This represents a significant reduction in the computational effort required for solving the problem (especially in 3D) and allows to consistently work with the full set of 12 FCC slip systems.
Numerical examples show that damage-dependent micro-flexible boundary conditions are, opposed to simpler micro-free and micro-hard boundary conditions, able to capture a realistic hardening response as well as crack development and propagation across grain boundaries. 
It is also shown that the model can be used in order to simulate coalescence of pre-existing voids. 
Furthermore, fracture of polycrystals is inherently a three-dimensional problem and although the necessary mechanisms can be comprehensively demonstrated on two-dimensional examples, crack fronts develop differently in three-dimensional grain structures. We show that our model is able to produce three-dimensional crack fronts by applying it to a three-dimensional polycrystal.

The presented modeling framework opens up for future investigation of damage initiation and micro-crack formation under consideration of the influence of grain boundary behavior. 
The role of grain boundaries as obstacles for the transport of plastic slip and how to model varying degrees of resistance against slip transmission has been a central topic of the presented work. We have employed the simplification that grain boundaries act as slip sinks. In reality, however, the process of slip transmission from one grain to the other is much more complex and depends among other things on the misalignment of crystal orientations and the type of dislocations. In order to account for the process of (partial) slip transmission between the grains on each side of the grain boundary, it is possible to extend the model by methods that allow to account for jumps in the field equations such as cohesive elements.
In order to capture the plastic slip resistance of each slip system based on misalignment angles, the micro-flexible boundary condition could be derived for a model that accounts for the plastic slip gradient in each slip system separately, such as the one presented in \cite{Auth2024}. The authors want to point out that such models, especially in combination with complex per-slip system boundary conditions come at a largely increased computation cost, as well as a significantly larger set of parameters.

Finally, the presented model is built in order to capture damage initiation and formation of micro-cracks. An important mechanism during damage initiation is void nucleation, which often occurs at hard particles. Upon combining this work with a contact model, it should be possible to simulate void nucleation by debonding at hard particles, as for example observed by \citeauthor{Achouri2013} \cite{Achouri2013} and different localization patterns during void nucleation as discussed by \citeauthor{Noell2018} \cite{Noell2018}.
In order to simulate void coalescence leading up to component failure, it is often necessary to develop solvers that can trace unstable crack growth, such as arc-length solvers. Recently, dissipation based arc-length solvers have shown good performance for phase-field models and seem to be a promising approach for thermodynamically consistent models such as the one presented in this work \cite{Borjesson2022}, \cite{Bharali2022}.

\section*{Acknowledgements}
The work in this paper has been funded by the Swedish Research Council (Vetenskapsrådet) under the grant number 2018-04318.

\FloatBarrier
\newpage
\printbibliography

@article{Azman2022,
    title = {{4D characterisation of void nucleation, void growth and void coalescence using advanced void tracking algorithm on in situ X-ray tomographic data}},
    year = {2022},
    journal = {Materials Today Communications},
    author = {Azman, M. A. and Le Bourlot, C. and King, A. and Fabr{\`{e}}gue, D. and Maire, E.},
    pages = {103892},
    volume = {32},
    publisher = {Elsevier Ltd},
    url = {https://doi.org/10.1016/j.mtcomm.2022.103892},
    doi = {10.1016/j.mtcomm.2022.103892},
    issn = {23524928}
}

@article{Carlsson2017,
    title = {{A comparison of the primal and semi-dual variational formats of gradient-extended crystal inelasticity}},
    year = {2017},
    journal = {Computational Mechanics},
    author = {Carlsson, Kristoffer and Runesson, Kenneth and Larsson, Fredrik and Ekh, Magnus},
    number = {4},
    pages = {531--548},
    volume = {60},
    publisher = {Springer Berlin Heidelberg},
    doi = {10.1007/s00466-017-1419-y},
    issn = {01787675}
}

@article{HernandezPadilla2014,
    title = {{A coupled phase-field model for ductile fracture in crystal plasticity}},
    year = {2014},
    journal = {Pamm},
    author = {Hernandez Padilla, Carlos Alberto and Markert, Bernd},
    number = {AUGUST},
    pages = {441--442},
    volume = {14},
    url = {http://doi.wiley.com/10.1002/pamm.201410208},
    isbn = {1022310224},
    doi = {10.1002/pamm.201410208},
    issn = {16177061}
}

@article{Dahlberg2013,
    title = {{A deformation mechanism map for polycrystals modeled using strain gradient plasticity and interfaces that slide and separate}},
    year = {2013},
    journal = {International Journal of Plasticity},
    author = {Dahlberg, Carl F.O. and Faleskog, Jonas and Niordson, Christian F. and Legarth, Brian Nyvang},
    pages = {177--195},
    volume = {43},
    url = {http://dx.doi.org/10.1016/j.ijplas.2012.11.010},
    doi = {10.1016/j.ijplas.2012.11.010},
    issn = {07496419}
}

@article{Simo1988,
    title = {{A framework for finite strain elastoplasticity based on maximum plastic dissipation and the multiplicative decomposition: Part 1. Continuum formulation}},
    year = {1988},
    journal = {Computer Methods in Applied Mechanics and Engineering},
    author = {Simo, J.C.},
    number = {2},
    pages = {199--219},
    volume = {66}
}

@article{Borjesson2022,
    title = {{A generalised path-following solver for robust analysis of material failure}},
    year = {2022},
    journal = {Computational Mechanics},
    author = {B{\"{o}}rjesson, Elias and Remmers, Joris J.C. and Fagerstr{\"{o}}m, Martin},
    number = {2},
    pages = {437--450},
    volume = {70},
    publisher = {Springer Berlin Heidelberg},
    doi = {10.1007/s00466-022-02175-w},
    issn = {14320924}
}

@article{Wulfinghoff2013,
    title = {{A gradient plasticity grain boundary yield theory}},
    year = {2013},
    journal = {International Journal of Plasticity},
    author = {Wulfinghoff, Stephan and Bayerschen, Eric and B{\"{o}}hlke, Thomas},
    pages = {33--46},
    volume = {51},
    publisher = {Elsevier Ltd},
    url = {http://dx.doi.org/10.1016/j.ijplas.2013.07.001},
    doi = {10.1016/j.ijplas.2013.07.001}
}

@article{Gurtin2002,
    title = {{A gradient theory of single-crystal viscoplasticity that accounts for geometrically necessary dislocations}},
    year = {2002},
    journal = {Journal of the Mechanics and Physics of Solids},
    author = {Gurtin, Morton E.},
    number = {1},
    pages = {5--32},
    volume = {50},
    doi = {10.1016/S0022-5096(01)00104-1},
    issn = {00225096}
}

@article{Bharali2023,
    title = {{A micromorphic phase-field model for brittle and quasi-brittle fracture}},
    year = {2023},
    journal = {Computational Mechanics},
    author = {Bharali, Ritukesh and Larsson, Fredrik and J{\"{a}}nicke, Ralf},
    publisher = {Springer Berlin Heidelberg},
    doi = {10.1007/s00466-023-02380-1},
    issn = {14320924},
    arxivId = {2206.11583}
}

@article{Ekh2004,
    title = {{A model framework for anisotropic damage coupled to crystal (visco)plasticity}},
    year = {2004},
    journal = {International Journal of Plasticity},
    author = {Ekh, Magnus and Lillbacka, Robert and Runesson, Kenneth},
    number = {12},
    pages = {2143--2159},
    volume = {20},
    doi = {10.1016/j.ijplas.2004.04.007},
    issn = {07496419}
}

@article{Miehe2010_historyvariable,
    title = {{A phase field model for rate-independent crack propagation: Robust algorithmic implementation based on operator splits}},
    year = {2010},
    journal = {Computer Methods in Applied Mechanics and Engineering},
    author = {Miehe, Christian and Hofacker, Martina and Welschinger, Fabian},
    number = {45-48},
    pages = {2765--2778},
    volume = {199},
    publisher = {Elsevier B.V.},
    url = {http://dx.doi.org/10.1016/j.cma.2010.04.011},
    doi = {10.1016/j.cma.2010.04.011},
    issn = {00457825}
}

@article{Bharali2022,
    title = {{A robust monolithic solver for phase-field fracture integrated with fracture energy based arc-length method and under-relaxation}},
    year = {2022},
    journal = {Computer Methods in Applied Mechanics and Engineering},
    author = {Bharali, Ritukesh and Goswami, Somdatta and Anitescu, Cosmin and Rabczuk, Timon},
    pages = {114927},
    volume = {394},
    publisher = {Elsevier B.V.},
    url = {https://doi.org/10.1016/j.cma.2022.114927},
    doi = {10.1016/j.cma.2022.114927},
    issn = {00457825}
}

@article{Gurtin2008,
    title = {{A theory of grain boundaries that accounts automatically for grain misorientation and grain-boundary orientation}},
    year = {2008},
    journal = {Journal of the Mechanics and Physics of Solids},
    author = {Gurtin, Morton E.},
    pages = {640--662},
    volume = {56},
    doi = {10.1016/j.jmps.2007.05.002}
}

@article{Auth2024,
    title = {{A thermodynamic framework for ductile phase-field fracture and gradient-enhanced crystal plasticity}},
    year = {2024},
    journal = {European Journal of Mechanics / A Solids},
    author = {Auth, Kim Louisa and Brouzoulis, Jim and Ekh, Magnus},
    pages = {105418},
    volume = {108},
    publisher = {Elsevier Masson SAS},
    url = {https://doi.org/10.1016/j.euromechsol.2024.105418},
    doi = {10.1016/j.euromechsol.2024.105418},
    issn = {0997-7538}
}

@article{Ambrosio1990,
    title = {{Approximation of functional depending on jumps by elliptic functional via Gamma‐convergence}},
    year = {1990},
    journal = {Communications on Pure and Applied Mathematics},
    author = {Ambrosio, Luigi and Tortorelli, Vincenzo Maria},
    number = {8},
    pages = {999--1036},
    volume = {43},
    doi = {10.1002/cpa.3160430805},
    issn = {10970312}
}

@incollection{Alessi2018,
    title = {{Comparison of phase-field models of fracture coupled with plasticity}},
    year = {2018},
    booktitle = {Computational Methods in Applied Sciences},
    author = {Alessi, R. and Ambati, M. and Gerasimov, T. and Vidoli, S. and De Lorenzis, L.},
    pages = {1--21},
    volume = {46},
    isbn = {9783319608853},
    doi = {10.1007/978-3-319-60885-3{\_}1},
    issn = {18713033}
}

@article{McBride2016,
    title = {{Computational and theoretical aspects of a grain-boundary model at finite deformations}},
    year = {2016},
    journal = {Technische Mechanik},
    author = {McBride, A. T. and Gottschalk, D. and Reddy, B. D. and Wriggers, P. and Javili, A.},
    number = {1-2},
    pages = {102--119},
    volume = {36},
    doi = {10.24352/UB.OVGU-2017-013},
    issn = {02323869}
}

@article{Maloth2023,
    title = {{Coupled Crystal Plasticity Phase-Field Model for Ductile Fracture in Polycrystalline Microstructures}},
    year = {2023},
    journal = {International Journal for Multiscale Computational Engineering},
    author = {Maloth, Thirupathi and Ghosh, Somnath},
    number = {2},
    pages = {1--19},
    volume = {21},
    doi = {10.1615/IntJMultCompEng.2022042164},
    issn = {15431649}
}

@article{Kacher2014,
    title = {{Dislocation interactions with grain boundaries}},
    year = {2014},
    journal = {Current Opinion in Solid State and Materials Science},
    author = {Kacher, Josh and Eftink, B. P. and Cui, B. and Robertson, I. M.},
    number = {4},
    pages = {227--243},
    volume = {18},
    publisher = {Elsevier Ltd},
    url = {http://dx.doi.org/10.1016/j.cossms.2014.05.004},
    doi = {10.1016/j.cossms.2014.05.004},
    issn = {13590286}
}

@article{Noell2017,
    title = {{Do voids nucleate at grain boundaries during ductile rupture?}},
    year = {2017},
    journal = {Acta Materialia},
    author = {Noell, Philip and Carroll, Jay and Hattar, Khalid and Clark, Blythe and Boyce, Brad},
    pages = {103--114},
    volume = {137},
    publisher = {Acta Materialia Inc.},
    url = {http://dx.doi.org/10.1016/j.actamat.2017.07.004},
    doi = {10.1016/j.actamat.2017.07.004},
    issn = {13596454}
}

@article{Achouri2013,
    title = {{Experimental characterization and numerical modeling of micromechanical damage under different stress states}},
    year = {2013},
    journal = {Materials and Design},
    author = {Achouri, Mohamed and Germain, Guenael and Dal Santo, Philippe and Saidane, Delphine},
    pages = {207--222},
    volume = {50},
    publisher = {Elsevier Ltd},
    doi = {10.1016/j.matdes.2013.02.075},
    issn = {18734197}
}

@article{Pineau2016,
    title = {{Failure of metals I: Brittle and ductile fracture}},
    year = {2016},
    journal = {Acta Materialia},
    author = {Pineau, A. and Benzerga, A. A. and Pardoen, T.},
    pages = {424--483},
    volume = {107},
    publisher = {Acta Materialia Inc.},
    url = {http://dx.doi.org/10.1016/j.actamat.2015.12.034},
    doi = {10.1016/j.actamat.2015.12.034},
    issn = {13596454}
}

@misc{Ferrite.jl,
    title = {{Ferrite.jl - Finite element toolbox for Julia}},
    year = {2021},
    author = {Carlsson, Kristoffer and Ekre, Fredrik and {Contributors}},
    url = {https://github.com/Ferrite-FEM/Ferrite.jl}
}

@article{gmsh,
    title = {{Gmsh: A 3-D finite element mesh generator with built-in pre- and post-processing facilities}},
    year = {2009},
    journal = {International Journal for Numerical Methods in Engineering},
    author = {Geuzaine, Christophe and Remacle, Jean François},
    number = {11},
    pages = {1309--1331},
    volume = {79},
    doi = {10.1002/nme.2579},
    issn = {00295981}
}

@article{Ekh2007,
    title = {{Gradient crystal plasticity as part of the computational modelling of polycrystals}},
    year = {2007},
    journal = {International Journal for Numerical Methods in Engineering},
    author = {Ekh, Magnus and Grymer, M. and Runesson, K. and Svedberg, T.},
    number = {2},
    pages = {197--220},
    volume = {72},
    url = {https://doi.org/10.1002/nme.2015},
    doi = {10.1002/nme.2015}
}

@article{Ekh2011,
    title = {{Influence of grain boundary conditions on modeling of size-dependence in polycrystals}},
    year = {2011},
    journal = {Acta Mechanica},
    author = {Ekh, Magnus and Bargmann, Swantje and Grymer, Mikkel},
    number = {1-2},
    pages = {103--113},
    volume = {218},
    doi = {10.1007/s00707-010-0403-9},
    issn = {00015970}
}

@article{Yalcinkaya2019,
    title = {{Inter-granular cracking through strain gradient crystal plasticity and cohesive zone modeling approaches}},
    year = {2019},
    journal = {Theoretical and Applied Fracture Mechanics},
    author = {Yal{\c{c}}inkaya, T and {\"{O}}zdemir, İ and Firat, A O},
    number = {July},
    pages = {102306},
    volume = {103},
    publisher = {Elsevier},
    url = {https://doi.org/10.1016/j.tafmec.2019.102306},
    doi = {10.1016/j.tafmec.2019.102306},
    issn = {0167-8442}
}

@article{Bezanson2017,
    title = {{Julia: A fresh approach to numerical computing}},
    year = {2017},
    journal = {SIAM Review},
    author = {Bezanson, Jeff and Edelman, Alan and Karpinski, Stefan and Shah, Viral B.},
    number = {1},
    pages = {65--98},
    volume = {59},
    doi = {10.1137/141000671},
    issn = {00361445},
    arxivId = {1411.1607}
}

@article{Quey2011,
    title = {{Large-scale 3D random polycrystals for the finite element method: Generation, meshing and remeshing}},
    year = {2011},
    journal = {Computer Methods in Applied Mechanics and Engineering},
    author = {Quey, R. and Dawson, P. R. and Barbe, F.},
    number = {17-20},
    pages = {1729--1745},
    volume = {200},
    publisher = {Elsevier B.V.},
    url = {http://dx.doi.org/10.1016/j.cma.2011.01.002},
    doi = {10.1016/j.cma.2011.01.002},
    issn = {00457825}
}

@article{Danisch2021,
    title = {{Makie.jl: Flexible high-performance data visualization for Julia}},
    year = {2021},
    journal = {Journal of Open Source Software},
    author = {Danisch, Simon and Krumbiegel, Julius},
    number = {65},
    pages = {3349},
    volume = {6},
    doi = {10.21105/joss.03349}
}

@article{Forest2009,
    title = {{Micromorphic Approach for Gradient Elasticity, Viscoplasticity, and Damage}},
    year = {2009},
    journal = {Journal of Engineering Mechanics},
    author = {Forest, Samuel},
    number = {3},
    pages = {117--131},
    volume = {135},
    doi = {10.1061/(asce)0733-9399(2009)135:3(117)},
    issn = {0733-9399}
}

@article{Aslan2011,
    title = {{Micromorphic approach to single crystal plasticity and damage}},
    year = {2011},
    journal = {International Journal of Engineering Science},
    author = {Aslan, O. and Cordero, N. M. and Gaubert, A. and Forest, S.},
    number = {12},
    pages = {1311--1325},
    volume = {49},
    publisher = {Elsevier Ltd},
    url = {http://dx.doi.org/10.1016/j.ijengsci.2011.03.008},
    doi = {10.1016/j.ijengsci.2011.03.008},
    issn = {00207225}
}

@article{Yalcinkaya2021,
    title = {{Misorientation and grain boundary orientation dependent grain boundary response in polycrystalline plasticity}},
    year = {2021},
    journal = {Computational Mechanics},
    author = {Yal{\c{c}}inkaya, Tuncay and {\"{O}}zdemir, İzzet and Tarik Tando{\u{g}}an, İzzet},
    number = {3},
    pages = {937--954},
    volume = {67},
    doi = {10.1007/s00466-021-01972-z},
    issn = {14320924}
}

@article{Bargmann2010,
    title = {{Modeling of polycrystals with gradient crystal plasticity: A comparison of strategies}},
    year = {2010},
    journal = {Philosophical Magazine},
    author = {Bargmann, Swantje and Ekh, Magnus and Runesson, Kenneth and Svendsen, Bob},
    number = {10},
    pages = {1263--1288},
    volume = {90},
    doi = {10.1080/14786430903334332},
    issn = {14786435}
}

@article{Svendsen2010,
    title = {{On the continuum thermodynamic rate variational formulation of models for extended crystal plasticity at large deformation}},
    year = {2010},
    journal = {Journal of the Mechanics and Physics of Solids},
    author = {Svendsen, Bob and Bargmann, Swantje},
    number = {9},
    pages = {1253--1271},
    volume = {58},
    publisher = {Elsevier},
    url = {http://dx.doi.org/10.1016/j.jmps.2010.06.005},
    doi = {10.1016/j.jmps.2010.06.005}
}

@article{Ambati2015,
    title = {{Phase-field modeling of ductile fracture}},
    year = {2015},
    journal = {Computational Mechanics},
    author = {Ambati, M. and Gerasimov, T. and De Lorenzis, L.},
    number = {5},
    pages = {1017--1040},
    volume = {55},
    publisher = {Springer Berlin Heidelberg},
    url = {http://dx.doi.org/10.1007/s00466-015-1151-4},
    doi = {10.1007/s00466-015-1151-4},
    issn = {01787675}
}

@article{Miehe2017,
    title = {{Phase-field modeling of ductile fracture at finite strains: A robust variational-based numerical implementation of a gradient-extended theory by micro}},
    year = {2017},
    journal = {International Journal for Numerical Methods in Engineering},
    author = {Miehe, Christian and Aldakheel, Fadi and Teichtmeister, Stephan},
    number = {9},
    pages = {816--863},
    volume = {111},
    url = {https://onlinelibrary.wiley.com/doi/10.1002/nme.5484},
    doi = {10.1002/nme.5484}
}

@article{DeLorenzis2016,
    title = {{Phase-field modelling of fracture in single crystal plasticity}},
    year = {2016},
    journal = {GAMM Mitteilungen},
    author = {De Lorenzis, L. and McBride, A. and Reddy, B. D.},
    number = {1},
    pages = {7--34},
    volume = {39},
    doi = {10.1002/gamm.201610002},
    issn = {09367195}
}

@article{Bayerschen2015,
    title = {{Review on slip transmission criteria in experiments and crystal plasticity models}},
    year = {2015},
    journal = {Journal of Materials Science},
    author = {Bayerschen, E. and McBride, a. T. and Reddy, B. D. and B{\"{o}}hlke, T.},
    url = {http://link.springer.com/10.1007/s10853-015-9553-4},
    doi = {10.1007/s10853-015-9553-4},
    issn = {0022-2461}
}

@article{Evers2004,
    title = {{Scale dependent crystal plasticity framework with dislocation density and grain boundary effects}},
    year = {2004},
    journal = {International Journal of Solids and Structures},
    author = {Evers, L. P. and Brekelmans, W. A.M. and Geers, M. G.D.},
    number = {18-19},
    pages = {5209--5230},
    volume = {41},
    doi = {10.1016/j.ijsolstr.2004.04.021},
    issn = {00207683}
}

@article{Fredriksson2005,
    title = {{Size-dependent yield strength of thin films}},
    year = {2005},
    journal = {International Journal of Plasticity},
    author = {Fredriksson, P. and Gudmundson, P.},
    number = {9},
    pages = {1834--1854},
    volume = {21},
    doi = {10.1016/j.ijplas.2004.09.005},
    issn = {07496419}
}

@article{Flouriot2003,
    title = {{Strain localization at the crack tip in single crystal CT specimens under monotonous loading: 3D Finite Element analyses and application to nickel-base superalloys}},
    year = {2003},
    journal = {International Journal of Fracture},
    author = {Flouriot, S. and Forest, S. and Cailletaud, G. and K{\"{o}}ster, A. and R{\'{e}}my, L. and Burgardt, B. and Gros, V. and Mosset, S. and Delautre, J.},
    number = {1-2},
    pages = {43--77},
    volume = {124},
    doi = {10.1023/B:FRAC.0000009300.70477.ba},
    issn = {03769429}
}

@article{Tensors.jl,
    title = {{Tensors.jl - Tensor computations in Julia}},
    year = {2019},
    journal = {Journal of Open Research Software},
    author = {Carlsson, Kristoffer and Ekre, Fredrik},
    number = {1},
    pages = {2--6},
    volume = {7},
    doi = {10.5334/jors.182},
    issn = {20499647}
}

@article{Noell2018,
    title = {{The mechanisms of ductile rupture}},
    year = {2018},
    journal = {Acta Materialia},
    author = {Noell, Philip J. and Carroll, Jay D. and Boyce, Brad L.},
    pages = {83--98},
    volume = {161},
    publisher = {Acta Materialia Inc.},
    url = {http://dx.doi.org/10.1016/j.actamat.2018.09.006},
    doi = {10.1016/j.actamat.2018.09.006},
    issn = {13596454}
}

@article{Aifantis2005,
    title = {{The role of interfaces in enhancing the yield strength of composites and polycrystals}},
    year = {2005},
    journal = {Journal of the Mechanics and Physics of Solids},
    author = {Aifantis, K. E. and Willis, J. R.},
    number = {5},
    pages = {1047--1070},
    volume = {53},
    doi = {10.1016/j.jmps.2004.12.003},
    issn = {00225096}
}

@article{Coleman1963,
    title = {{The thermodynamics of elastic materials with heat conduction and viscosity}},
    year = {1963},
    journal = {Archive for Rational Mechanics and Analysis},
    author = {Coleman, Bernard D and Noll, Walter},
    pages = {167--178},
    volume = {13},
    publisher = {Springer},
    doi = {10.1007/BF01262690}
}

@article{Rovinelli2018,
    title = {{Using machine learning and a data-driven approach to identify the small fatigue crack driving force in polycrystalline materials}},
    year = {2018},
    journal = {npj Computational Materials},
    author = {Rovinelli, Andrea and Sangid, Michael D. and Proudhon, Henry and Ludwig, Wolfgang},
    number = {1},
    pages = {1--10},
    volume = {4},
    publisher = {Springer US},
    url = {http://dx.doi.org/10.1038/s41524-018-0094-7},
    doi = {10.1038/s41524-018-0094-7},
    issn = {20573960}
}

@article{Spannraft2020GrainDecohesion,
    title = {{Grain boundary interaction based on gradient crystal inelasticity and decohesion}},
    year = {2020},
    journal = {Computational Materials Science},
    author = {Spannraft, Lucie and Ekh, Magnus and Larsson, Fredrik and Runesson, Kenneth and Steinmann, Paul},
    number = {October 2019},
    pages = {109604},
    volume = {178},
    publisher = {Elsevier},
    url = {https://doi.org/10.1016/j.commatsci.2020.109604},
    doi = {10.1016/j.commatsci.2020.109604},
    issn = {09270256}
}

\newpage
\FloatBarrier
\begin{appendices}
\section{Material parameters}
\label{section:appendix_parameters}
Material parameters for the two-dimensional simulations are given in Table \ref{tab:material_parameters}. The length scales for the three-dimensional simulation deviate such that $l_\mathrm{g}=0.1333\,L$ and $\ell_0=0.08\,L$. The other parameters are adjusted as presented in Table \ref{tab:material_parameters}. All simulations use full FCC slip-systems, see Table \ref{tab:slipsystems} for the unit cell slip systems. In order to account for the different crystal orientations in the grains, the unit cell of each grain is rotated according to the Rodrigues parameters given in Table \ref{tab:rodrigues_angles}.
\begin{table}[h!]
    \centering
{\linespread{1.2}\selectfont
    \begin{tabular}{||l l r l||}
        \hline
        Parameter &  & Value & Unit \\ [0.5ex]
        \hline \hline
        Bulk modulus & $\kappa$ & 71660 & MPa \\
        \hline
        Shear modulus & $\mu$ & 27260  & MPa \\ [0.5ex]
        \hline \hline
        Yield stress & $\tau^\mathrm{y}$ & 345 & MPa \\ [0.5ex]
        \hline
        Isotropic hardening modulus & $H_\mathrm{iso}$ & 250 & MPa \\
        \hline
        Gradient hardening modulus & $H^\mathrm{g}$ & 1000 & MPa \\
        \hline
        Gradient hardening length scale & $l_\mathrm{g}$ & $0.0533\,L$ &  \\
        \hline \hline
        Visco-plastic relaxation time & $t^\ast$ & 1 & s \\
        \hline
        Visco-plastic drag stress & $\sigma^\mathrm{d}$ & 500 & MPa \\
        \hline
        Visco-plastic exponent & $m$ & 8 & - \\ [0.5ex]
        \hline
        \hline
        Initial micro-flexiblity & $C_{\Gamma,0}$ & $1 \, / \left(H^\mathrm{g} \, l_\mathrm{g}^2\right)$  &  \\
        \hline
        Micro-flexibilty slope & $C_{\Gamma}^\mathrm{d}$ & $20 \,  / \left(H^\mathrm{g} \, l_\mathrm{g}^2\right)$  &  \\ [0.5ex]
        \hline
        \hline
        Effective fracture energy & $\mathcal{G}^\mathrm{d}_0 / \ell_0$ & 300 & $\mathrm{N}/\mathrm{mm}^2$ \\
        \hline
        Phase-field length scale & $\ell_0$ & $0.02\,L$ &  \\
        \hline
        Micromorphic penalty parameter & $\alpha$ & $200\, \mathcal{G}^\mathrm{d}_0\,/\,\ell_0$ & \\
        \hline
        Critical plastic strain & $\epsilon^\mathrm{p}_\mathrm{crit}$ & 0.1 & - \\
        \hline
        Degradation exponent & $n$ & 2 & - \\
        \hline
    \end{tabular}
    }
    \caption{Base material parameters employed for the numerical experiments. $L$ is a global length scale for the simulated structure and results are independent of $L$ upon scaling. Size independence of the model can be derived in the same way as demonstrated in \cite{Auth2024}.}
    \label{tab:material_parameters}
\end{table}
\begin{table}[h!]
    \centering
{\linespread{1.2}\selectfont
    \begin{tabular}{||c c c|c c c|c c c|c c c||}
        \hline
       $\alpha$  & $\bar{\boldsymbol{s}}_\alpha$ & $\bar{\boldsymbol{m}}_\alpha$ &
       $\alpha$  & $\bar{\boldsymbol{s}}_\alpha$ & $\bar{\boldsymbol{m}}_\alpha$ & 
       $\alpha$  & $\bar{\boldsymbol{s}}_\alpha$ & $\bar{\boldsymbol{m}}_\alpha$ &
       $\alpha$  & $\bar{\boldsymbol{s}}_\alpha$ & $\bar{\boldsymbol{m}}_\alpha$ \\ [0.5ex]
       \hline \hline
       1 & $\left[\bar{1} \, 1 \, 0 \right]$ & $\left[1 \, 1 \, 1 \right]$ &
       4 & $\left[\bar{1} \, \bar{1} \, 0 \right]$ & $\left[1 \, \bar{1} \, \bar{1} \right]$ &
       7 & $\left[1 \, 1 \, 0 \right]$ & $\left[\bar{1} \, 1 \, \bar{1} \right]$ &
       10 & $\left[1 \, \bar{1} \, 0 \right]$ & $\left[\bar{1} \, \bar{1} \, 1 \right]$ \\
       \hline
       2 & $\left[1 \, 0 \, \bar{1} \right]$ & $\left[1 \, 1 \, 1 \right]$ &
       5 & $\left[1 \, 0 \, 1 \right]$ & $\left[1 \, \bar{1} \, \bar{1} \right]$ &
       8 & $\left[\bar{1} \, 0 \, 1 \right]$ & $\left[\bar{1} \, 1 \, \bar{1} \right]$ &
       11 & $\left[\bar{1} \, 0 \, \bar{1} \right]$ & $\left[\bar{1} \, \bar{1} \, 1 \right]$ \\
       \hline
       3 & $\left[0 \, \bar{1} \, 1 \right]$ & $\left[1 \, 1 \, 1 \right]$ &
       6 & $\left[0 \, 1 \, \bar{1} \right]$ & $\left[1 \, \bar{1} \, \bar{1} \right]$ &
       9 & $\left[0 \, \bar{1} \, \bar{1} \right]$ & $\left[\bar{1} \, 1 \, \bar{1} \right]$ &
       12 & $\left[0 \, 1 \, 1 \right]$ & $\left[\bar{1} \, \bar{1} \, 1 \right]$ \\
       \hline
    \end{tabular}
}
    \caption{FCC slip systems. The unit cell of each grain in the numerical examples is rotated as described in Table \ref{tab:rodrigues_angles}.}
    \label{tab:slipsystems}
\end{table}

\begin{table}[htb]
    \centering
{\linespread{1.2}\selectfont
    \begin{tabular}{||c|c c c||c|c c c||}
        \hline
        Grain & \multicolumn{3}{c||}{2D Rodrigues angles} & Grain & \multicolumn{3}{c||}{2D Rodrigues angles}\\ [0.5ex]
        \hline \hline
        1 & -3.67 & 2.54 & -0.52 & 6 & 1.12 & -1.05 & -1.1 \\
        \hline
        2 & 3.8 & 22.4 & 16.75 & 7 & 1.4 & 2.07 & -0.64 \\
        \hline
        3 & -6.32 & 5.63 & -27.73 & 8 & 0.55 & -0.53 & -0.55 \\
        \hline
        4 & -0.22 & 0.26 & -0.1 & 9 & 0.71 & -0.94 & 1.02 \\
        \hline
        5 & 0.08 & 0.54 & 0.02 & 10 & -0.97 & -0.42 & 0.37 \\
        \hline
    \end{tabular}
}
    \quad
{\linespread{1.2}\selectfont
    \begin{tabular}{||c|c c c||}
    \hline
    Grain & \multicolumn{3}{c||}{3D Rodrigues angles} \\ [0.5ex]
    \hline \hline
    1 &0.70  &  1.81 &  0.72 \\
    \hline
    2 &  -1.76  & 1.73  & 4.43 \\
    \hline
    3 &  2.56 & -1.54 & -0.23 \\ 
    \hline
    4 & -0.35 &  0.09 & -0.32 \\
    \hline
    \end{tabular}
}
    \caption{Rodrigues angles for the rotation of the FCC unit cell in each grain. The left table shows the angles for the 2D grain structure with 6 grains, the right table shows the angles for the 3D grain structure with 4 grains.}
    \label{tab:rodrigues_angles}
\end{table}

The loading rate at the upper side of the structures is $0.05\,L \; / \; t^\ast$. The sides are loaded accordingly, in order to achieve displacement controlled simple shear loading conditions.
\section{Geometry and Meshing details}
\label{appendix:geometry_meshing}
The polycrystals are generated by the open-source software 
\texttt{Neper} \cite{Quey2011}. The prompts for generating the two- and three-dimensional tessalations are \texttt{neper -T -dim 2 -n 10 -id 2} and \texttt{neper -T -dim 3 -n 4 -id 1}, respectively.
The voids in the two-dimensional structure as placed at positions $\left[0.36\,L, 0.415\,L\right]$ and $\left[0.49\,L, 0.425\,L\right]$ and have a diameter of $0.045\,L$.
All meshes are generated by \texttt{gmsh} \cite{gmsh}. The average element sizes in the different areas of the meshes are given in Table \ref{table:element_sizes}.
\begin{table}[htb]
    \centering
{\linespread{1.2}\selectfont
    \begin{tabular}{|| l | c c c ||}
        \hline
        Element size & 2D polycrystal & 2D polycrystal with voids & 3D polycrystal \\ [0.5ex]
         \hline \hline
        Background mesh & $0.0125\,L$ & $0.0125\,L$ & $0.053333\,L$ \\
        \hline
        Grain boundaries &$0.00875\,L$ & $0.00875\,L$ & $0.040000\,L$\\
        \hline
        Grain boundary intersections & $0.005\,L$ & $0.005\,L$ & $0.040000\,L$\\
        \hline
        Voids & - & $0.0005625\,L$ & - \\
        \hline
        Between voids & - & $0.001125\,L$ & - \\
        \hline
        Grain boundary right of voids & - & $0.00225\,L$ & - \\ [0.5ex]
        \hline \hline
        Number of nodes & $12\,146$ & $24\,033$ & $16\,065$ \\
        \hline
        Number of Triangles/Tetrahedra & $22\,616$ & $45\,688$ & $76\,506$ \\
        \hline
    \end{tabular}
}
    \caption{Mesh characteristics: Average element sizes in the different areas of the meshes and number of nodes/elements.}
    \label{table:element_sizes}
\end{table}

\section{Time integration}
\label{appendix:time_integration}
By applying backward Euler time integration to Equations (\ref{eq:evolution_Fp}) and (\ref{eq:evolution_k}) and discretizing Equations (\ref{eq:evolution_ep}) and (\ref{eq:varphi_KKT}), time-discretized forms of the local variables are obtained
\begin{align}
    \,^{n+1}\boldsymbol{F}_\mathrm{p}^{-1}
    &=
    \,^{n}\boldsymbol{F}_\mathrm{p}^{-1}
    \cdot
    \left(
        \boldsymbol{I} - \sumalpha \frac{\,^{n+1}\Delta\lambda_\alpha}{g_\mathrm{e}\left(\,^\mathrm{n+1}\varphi,\,^{n+1}\epsilon^\mathrm{p}\right)} \, \left(\bar{\boldsymbol{s}}_\alpha \otimes \bar{\boldsymbol{m}}_\alpha\right) \, \mathrm{sign}\left(\,^{n+1}\tau_\alpha\right)
    \right)
    \\
    \,^{n+1}k_\alpha
    &=
    \,^{n}k_\alpha - \,^{n+1}\Delta \lambda_\alpha
    \\
    \,^{n+1}\epsilon^\mathrm{p}
    &=
    \,^{n}\epsilon^\mathrm{p} + \sqrt{\sumalpha \,^{n+1}\Delta \lambda_\alpha^2}
    \\
    \label{eq:n+1varphi}
    \,^{n+1}\varphi
    &=
    \mathrm{max}\left(\,^{n}\varphi, \, \,^{n+1}\varphi^\mathrm{trial}\right) \,.
\end{align}
The signs of the Schmid stresses $\,^{n+1}\tau_\alpha$ can be computed based on elastic trial stresses. Since all local variables can then be expressed in terms of $\,^{n+1}\Delta \lambda_\alpha$ and $\,^{n+1}\varphi^\mathrm{trial}$, it remains the solve the coupled equation system resulting from Equations (\ref{eq:lambda_alpha}) and (\ref{eq:local_residual_pf})
\begin{align}
    \label{eq:residual_Deltalambda}
    \mathcal{R}_{\Delta \lambda_\alpha}\left(\,^{n+1}\Delta\lambda_\alpha,\,\,^{n+1}\varphi^\mathrm{trial}\right)
    &=
    \,^{n+1}\Delta\lambda_\alpha - \frac{\Delta t}{t^\ast} \,
    \left<
        \frac{
            \Phi_\alpha \left(
                \,^{n+1}\Delta\lambda_\alpha,\,\,^{n+1}\varphi^\mathrm{trial}
            \right)
        }{
            \sigma_\mathrm{d}
        }
    \right>^m = 0
    \\
    \label{eq:residual_varphi}
    \mathcal{R}_{\varphi^\mathrm{trial}}\left(\,^{n+1}\Delta\lambda_\alpha,\,\,^{n+1}\varphi^\mathrm{trial}\right)
    &=
    - \frac{\partial g_\mathrm{e}\left(
        \,^{n+1}\varphi^\mathrm{trial} ,\,\,^{n+1}\Delta\lambda_\alpha
    \right)
    }{\partial \varphi}
    \hat{\Psi}_\mathrm{e}^{+}
    -
    \frac{\mathcal{G}_0^\mathrm{d}}{\ell_0} \,\,^{n+1}\varphi^\mathrm{trial} 
    -
    \alpha \, \left(\,^{n+1}\varphi^\mathrm{trial}-d\right) = 0
    \,.
\end{align}
Notice that while $\mathcal{R}_{\Delta\lambda_\alpha}$ implicitly depends on $\,^{n+1}\varphi^\mathrm{trial}$ via Equation (\ref{eq:n+1varphi}), $\mathcal{R}_{\varphi^\mathrm{trial}}$ explicitly depends on it, meaning it solves $Y_\varphi(\varphi^\mathrm{trial})=0$.

\section{Staggered solver}
\label{appendix:solver}
\begin{figure}[htb]
    \centering
    \includegraphics[]{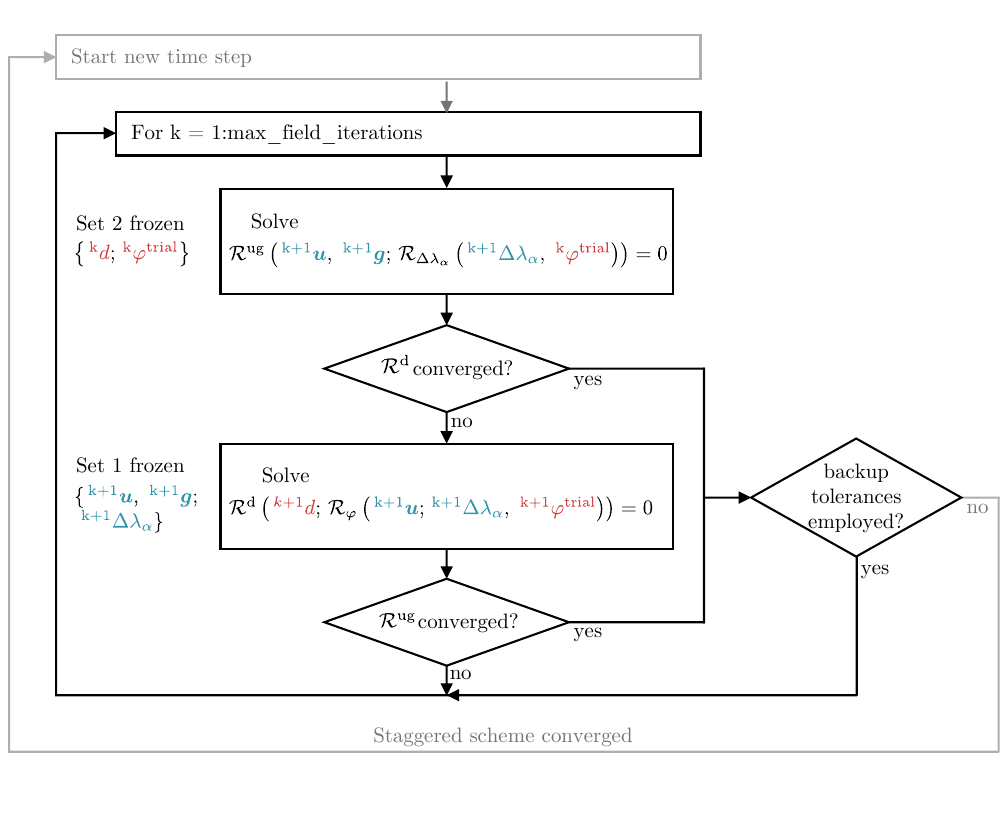}
    \caption{Staggered solver: The residual equations are solved in two steps, where the first step ($\mathcal{R}^\mathrm{ug}$) consists of the equilibrium equation and the gradient plasticity equation and the second step consists of the phase-field equation ($\mathcal{R}^\mathrm{d}$). Iterations between the two steps go back and forth until all global residual equations fulfill the final convergence criterion \textit{for the same global variables}. During the iteration process, unnecessary iterations are cut by allowing each single step to converge according to back-up tolerances. However, convergence according to the back-up tolerances does not allow to leave the field iteration loop.}
    \label{fig:staggered_scheme}
\end{figure}

\begin{table}[htb]
    \centering
{\linespread{1.2}\selectfont
    \begin{tabular}{||c | l c l l||}
        \hline
         ~ & ~ & Weak form & \vtop{\hbox{\strut Residual}\hbox{\strut tolerance}}& \vtop{\hbox{\strut Field }\hbox{\strut tolerance}} \\ [0.5ex]
        \hline \hline
        \multirow{3}{4em}{Final tolerances} &
        Displacement field $\boldsymbol{u}$ & Equation \ref{eq:weak_equil} & $L^2 \cdot 10^{-6}~\mathrm{N mm^{-2}} $ & $L \cdot 10^{-8}$ \\
        \cline{2-5}
        & Hardening strain gradient $\boldsymbol{g}$ & Equation \ref{eq:weak_galpha} &
        $L^2 \cdot 10^{-11} $ & $L^{-1} \cdot 10^{-6}$ \\
        \cline{2-5}
        & Global phase-field $d$ & Equation \ref{eq:weak_pf} &
        $L^3 \cdot 10^{-10}~\mathrm{N mm^{-3}} $ & $10^{-8}$ \\ [0.5ex]
        \hline \hline
        \multirow{3}{4em}{Back-up tolerances} &
        Displacement field $\boldsymbol{u}$ & Equation \ref{eq:weak_equil} & $L^2 \cdot 10^{-3}~\mathrm{N mm^{-2}} $ & $L \cdot 10^{-8}$ \\
        \cline{2-5}
        & Hardening strain gradient $\boldsymbol{g}$ & Equation \ref{eq:weak_galpha} &
        $L^2 \cdot 10^{-8} $ & $L^{-1} \cdot 10^{-6}$ \\
        \cline{2-5}
        & Global phase-field $d$ & Equation \ref{eq:weak_pf} &
        $L^3 \cdot 10^{-10}~\mathrm{N mm^{-3}} $ & $10^{-8}$
        \\
        \hline
    \end{tabular}
}
    \caption{Field-wise tolerances applied for global convergence examination. A field is considered converged when either of the residual or the field tolerance (relating to the update of field variables during a single Newton-correction) is fulfilled. $L$ is assumed to have unit mm.}
    \label{tab:global_tolerances}
\end{table}

A staggered solver is employed to solve the coupled equation system. Equations (\ref{eq:weak_equil}) and (\ref{eq:weak_galpha}) are solved in a monolithic way in the first step, referred to by $\mathcal{R}^\mathrm{ug}$ in Figure \ref{fig:staggered_scheme}, Equation (\ref{eq:weak_pf}) is solved in the second step, referred to by $\mathcal{R^\mathrm{d}}$, in each staggered iteration. We apply field-wise tolerances and each field can converge either based on a tolerance on its residual equation (residual tolerance) or on the update of the variable vector within a single Newton iteration (field tolerance).
The solver algorithm is represented in Figure \ref{fig:staggered_scheme}. In the same manner as presented in \cite{Auth2024}, the staggered scheme extends to the local level. Therefore, the global fields $\boldsymbol{u}$, $\boldsymbol{g}$ and $d$, as well as the local variables $\Delta\lambda_\alpha$ and $\varphi^\mathrm{trial}$ are split into two sets. Set 1 includes $\{\boldsymbol{u},\,\boldsymbol{g},\,\Delta\lambda_\alpha\}$, which are the variables relating to $\mathcal{R}^\mathrm{ug}$ and Set 2 includes $\{d,\,\varphi^\mathrm{trial}\}$, relating to $\mathcal{R}^\mathrm{d}$. In each step only the set of variables relating to the solved equations can be updated, while the other set of variables is frozen to the values of the previous iteration. Since the two local variables are part of different sets, this means that each local residual equation is solved separately and only in one of the two staggered steps.

In each staggered step, the global and local residual equations are solved by Newton iterations. After solving the system corresponding to the respective staggered step, convergence of the global residual equation corresponding to the other staggered step is checked. If both global residual equations fulfill the convergence criteria without further updates of the variables, the time step is converged.

\begin{table}[htb]
    \centering
{\linespread{1.2}\selectfont
    \begin{tabular}{|| l | r ||}
        \hline
        Solver parameter & value \\ [0.5ex]
        \hline \hline
        \texttt{max\_field\_iter} & 500 \\
        \hline
        \texttt{max\_iter} & 25 \\
        \hline
        \texttt{iter\_backuptol} & 8 \\
        \hline
        \texttt{n\_iter\_backuptols} for $\mathcal{R}^\mathrm{ug} $ & 25 \\
        \hline
        \texttt{n\_iter\_backuptols} for $\mathcal{R}^\mathrm{d} $ & 2 \\
        \hline
        \texttt{max\_field\_div} & 2 \\
        \hline
        \texttt{max\_divergence\_count} & 5 \\
        \hline
        \texttt{max\_time\_refinement\_level} & 20 \\
        \hline
        \texttt{n\_field\_iter\_coarsen} & 20 \\
        \hline
    \end{tabular}
}
    \caption{Parameters employed in the staggered solver.}
    \label{tab:solver_params}
\end{table}

The convergence criteria for the global residual equations are given in Table \ref{tab:global_tolerances}.
The iterations on the local residual equation employ a fixed tolerance of $10^{-8}$. 
In order to cut unnecessary iterations when solving each of the global residual equations, coarser back-up tolerances can be used (early) in the field iterations. However, the use of the back-up tolerances always triggers further iterations between the fields and thus a time step can never finally converge based on the back-up tolerances. The use of the back-up tolerances can be triggered by any of the following scenarios:
\begin{itemize}
    \item The residual tolerance is still larger than the back-up tolerance in a set iteration \texttt{iter\_backuptol}.
    \item Either of the previous two solves (one on each step) has employed more than \texttt{n\_iter\_backuptols} iterations AND the first computed residual is larger than the back-up tolerance. \texttt{n\_iter\_backuptols} is chosen for each step separately.
    \item The maximum number of iterations \texttt{max\_iter} is reached, but the final tolerances are not fulfilled.
\end{itemize}
These solver parameters are given in Table \ref{tab:solver_params}.
Additionally, iterations are terminated when divergence is recognized as follows:
\begin{itemize}
    \item Field iterations: The first residual in a solve is larger than the first residual of the previous solve more than \texttt{max\_field\_divergence} times in a row (for the same residual equation).
    \item Newton iterations for each solve: The residual has increased more than \texttt{max\_divergence\_count} subsequent times.
\end{itemize}
Within the Newton iterations, a basic line search is applied. Upon increasing residuals, the computed variable update $\Delta \boldsymbol{a}$ is reduced such that $\Delta \boldsymbol{a} = 0.7^\text{number of increases} \, \Delta \boldsymbol{a}_\text{original}$.

When a time step fails (e.g. due to recognized divergence or lack of material model convergence), the time step is halved. Time step halving is allowed up to \texttt{max\_time\_refinement\_level} times. 
After \texttt{n\_timesteps\_recoarsen} converged time steps in a row, the time step is doubled if the previous time step did not use more than \texttt{n\_field\_iter\_coarsen} field iterations. Time step doubling is not allowed to lead to coarser time steps than the initial time step $\Delta t_\mathrm{inital}$.
\end{appendices}

\end{document}